\documentclass[reqno, huge]{amsart}
\usepackage{latexsym}
\usepackage{amsmath}
\usepackage{amssymb}
\usepackage{amsthm}

\usepackage{CJK}
\usepackage{amscd}
\usepackage{graphicx}
\usepackage{xcolor}
\usepackage[colorlinks,citecolor=blue]{hyperref}
\usepackage{tikz}
\usepackage{tikz-cd}
\usepackage[a4paper,top=3cm,bottom=3cm,left=3cm,right=3cm]{geometry}
\usepackage{cite}

\newtheorem{theorem}{Theorem}[section]
\newtheorem{corollary}[theorem]{Corollary}
\newtheorem{lemma}[theorem]{Lemma}
\newtheorem{proposition}[theorem]{Proposition}
\theoremstyle{definition}
\newtheorem{definition}[theorem]{Definition}
\theoremstyle{remark}

\numberwithin{equation}{section}

\begin{document}

\title[]{A note on Bridgeland Moduli Spaces and Moduli Spaces of sheaves on $X_{14}$ and $Y_3$}

\subjclass[2010]{Primary 14F05; secondary 14J45, 14D20, 14D23}
\keywords{Derived Categories, Bridgeland moduli spaces, Kuznetsov components}\thanks{This note is part of the Undergraduate Mathematcial Research Project of the first author mentored by the second author.The second author is supported by ERC Consolidator grant WallCrossAG, no. 819864. }

\author{Zhiyu Liu}
\address{Department of Mathematics, Sichuan University, Chengdu, Sichuan Province 610064 P. R.
China}
\email{zhiyuliu@stu.scu.edu.cn}

\author{Shizhuo Zhang}
\address{School of Mathematics, The University of Edinburgh, JCMB Building, Kings Building, Edinburgh, EH9 3FD}
\email{Shizhuo.Zhang@ed.ac.uk}

\begin{abstract}
We study Bridgeland moduli spaces of semistable objects of $(-1)$-classes and $(-4)$-classes in the Kuznetsov components on index one prime Fano threefold $X_{4d+2}$ of degree $4d+2$ and index two prime Fano threefold $Y_d$ of degree $d$ for $d=3,4,5$. For every Serre-invariant stability condition on the Kuznetsov components, we show that the moduli spaces of stable objects of $(-1)$-classes on $X_{4d+2}$ and $Y_d$ are isomorphic. We show that moduli spaces of stable objects of $(-1)$-classes on $X_{14}$ are realized by Fano surface $\mathcal{C}(X)$ of conics, moduli spaces of semistable sheaves $M_X(2,1,6)$ and $M_X(2,-1,6)$ and the correspondent moduli spaces on cubic threefold $Y_3$ are realized by moduli spaces of stable vector bundles $M^b_Y(2,1,2)$ and $M^b_Y(2,-1,2)$. We show that moduli spaces of semistable objects of $(-4)$-classes on $Y_{d}$ are isomorphic to the moduli spaces of instanton sheaves $M^{inst}_Y$ when $d\neq 1,2$, and show that there're open immersions of $M^{inst}_Y$ into moduli spaces of semistable objects of $(-4)$-classes when $d=1,2$. Finally, when $d=3,4,5$ we show that these moduli spaces are all isomorphic to $M^{ss}_X(2,0,4)$.


\end{abstract}

\maketitle

\setcounter{tocdepth}{1}

\section{Introduction}

\subsection{Background and main results}

The notion of stability on a triangulated category was introduced by Bridgeland in \cite{T03}. It enables us to construct moduli spaces of semistable objects in a general triangulated category, which provides a powerful machinery to study the geometry of classical moduli spaces. One of the most striking progress is made in  \cite{bayer2014projectivity} and \cite{bayer2014mmp}, where the authors construct moduli space of Bridgeland stable objects on a K3 surface $S$ with respect to a stability condition $\sigma$ in the stability manifold $\mathrm{Stab}(S)$. They realized each birational model of a moduli space of stable sheaves on such a K3 surface $S$ as Bridgeland moduli space. They showed that the Minimal model program for this moduli space corresponds to a wall-crossing in $\mathrm{Stab}(S)$ in a precise way. On the other hand, it has been widely accepted that a non-trivial semi-orthogonal component $\mathcal{K}u(X)$ (called the Kuznetsov component) of a bounded derived category $\mathrm{D^b}(X)$ of a smooth projective variety $X$ encodes essential information of its birational geometry and classical moduli spaces on it. Recently, in \cite{BLMS}, the authors provide a criterion to induce a stability condition $\sigma$ on the  Kuznetsov components $\mathcal{K}u(X)$ of a series prime Fano threefolds $X$ from weak stability conditions on $\mathrm{D^b}(X)$.  Thus constructing Bridgeland moduli spaces $\mathcal{M}_{\sigma}(\mathcal{K}u(X),c)$ of $\sigma$-stable objects of class $c$ in $\mathcal{K}u(X)$ to study birational geometry of classical moduli spaces of semistable sheaves on $X$ becomes possible. In the present note, we make an attempt in this direction.



Let $Y_d$ be a prime Fano threefold of index 2 and degree $d$. We consider the semiorthogonal decomposition of $Y_d$ given by 
\[\mathrm{D^b}(Y_d)=\langle \mathcal{K}u(Y_d), \mathcal{O}_{Y_d}, \mathcal{O}_{Y_d}(1)\rangle\]
The numerical Grothendieck group $\mathcal{N}(\mathcal{K}u(Y_d))$ of $\mathcal{K}u(Y_d)$ is a rank 2 lattice spanned by 
\[v:=[I_L]=1-\frac{1}{d}H^2, w:=H-\frac{1}{2}H^2+(\frac{1}{6}-\frac{1}{d})H^3\]
where $I_L$ is an ideal sheaf of a line $L$ on $Y_d$. 

Let $X_{2g-2}$ be an index one prime Fano threefold of genus $g$ and degree $2g-2$. The semi-orthogonal decompositions of $X_{2g-2}$ are given in \cite{kuz09} and \cite{kuz06}.
\[\mathrm{D^b}(X_{2g-2})=\langle\mathcal{A}_{X_{2g-2}}, \mathcal{O}_{X_{2g-2}}, \mathcal{E}^{\vee}\rangle\]

The  numerical Grothendieck group $\mathcal{N}(\mathcal{A}_{X_{2g-2}})$ of $\mathcal{A}_{X_{2g-2}}$ is a rank 2 lattice spanned by
\[s:=[I_C]=1-\frac{1}{g-1}H^2, t:=H-(\frac{g}{2}+1)L-\frac{16-g}{12}P\]
where $I_C$ is an ideal sheaf of a conic $C$ on $X_{2g-2}$. 

We call $u\in \mathcal{N}(\mathcal{K}u(Y_d))$ (or $\mathcal{N}(X_{2g-2})$) a $(-r)$-class if $\chi(u,u)=-r$. As noted in \cite[Proposition 3.9]{kuz09}, the Chern character map identify the numerical Grothendieck group with the lattice generated by Chern characters of some sheaves, then we will use notations of numerical class and Chern character alternatively.

In \cite{kuz09} and \cite{KPS}, the author established equivalences of triangulated categories $\mathcal{K}u(Y_d)\cong\mathcal{A}_{X_{4d+2}}$ for pair some $(Y_d, X_{4d+2})\in Z_d\subset \mathcal{MF}^2_{d}\times \mathcal{MF}^1_{4d+2}$ and for $3\leq d\leq 5$ that map $1-L$ to $1-2L$, and prove the isomorphisms of Fano surface $\Sigma(Y_d)$ of lines and Fano surface $\mathcal{C}(X_{4d+2})$ via the equivalence. In the present note, we mainly study various classical moduli spaces on them from a modern point of view. We apply the techniques developed in \cite{BMMS},\cite{PY20},\cite{APR} and \cite{zhang2020bridgeland} to show the following results in this note. We denote by $\mathcal{M}_{\sigma}(\mathcal{A}_X, c)$ the moduli space of $\sigma$-stable objects of numerical class $c$ in $\mathcal{A}_X$, and $\mathcal{M}^{ss}_{\sigma}(\mathcal{A}_X, c)$ for S-equivalence classes of $\sigma$-semistable objects. 

\begin{theorem} 
Let $(Y,X)\in Z_d\subset \mathcal{MF}^2_{d}\times \mathcal{MF}^1_{4d+2}$ and $3\leq d\leq 5$. Let $\sigma', \sigma$ be two Serre-invariant stability conditions on $\mathcal{A}_X$ and $\mathcal{K}u(Y)$ respectively. Then the equivalence $\mathcal{K}u(Y)\cong \mathcal{A}_X$ induces following isomorphisms between Bridgeland moduli spaces:

\begin{enumerate}
    \item \emph{(Corollary \ref{all moduli spaces of (-1)-class isom})}  $\mathcal{M}_{\sigma}(\mathcal{K}u(Y), a) \xrightarrow{\cong} \mathcal{M}_{\sigma'}(\mathcal{A}_X, a')$ where  $a',a$ are any $(-1)$-classes in $\mathcal{A}_X$ and $\mathcal{K}u(Y)$ respectively.
    
    \item \emph{(Corollary \ref{all moduli spaces of (-4)-class isom})} When $d\neq 4$, we have  $\mathcal{M}^{ss}_{\sigma}(\mathcal{K}u(Y), b) \xrightarrow{\cong} \mathcal{M}^{ss}_{\sigma'}(\mathcal{A}_X, b')$ where $b',b$ are any two $(-4)$-classes in $\mathcal{A}_X$ and $\mathcal{K}u(Y)$ respectively. 
    
    \item \emph{(Theorem \ref{iso-4})} When $d=4$, we have  $\mathcal{M}^{ss}_{\sigma}(\mathcal{K}u(Y), 2-2L) \xrightarrow{\cong} \mathcal{M}^{ss}_{\sigma'}(\mathcal{A}_X, 2-4L)$.
\end{enumerate}
\end{theorem}

If $d=3$, up to sign, there are three $(-1)$-classes $v, v-w$ and $2v-w$ in $\mathcal{N}(\mathcal{K}u(Y_3))$ and another three $(-1)$-classes $s, 3s-t, 2s-t$ in $\mathcal{N}(\mathcal{A}_{X_{14}})$, respectively. All the Bridgeland moduli spaces of stable objects of $(-1)$-classes in $\mathcal{K}u(Y_3)$ and $\mathcal{A}_{X_{14}}$ can be realized as classical moduli spaces on $Y_3$ and $X_{14}$:

\begin{theorem}\label{negativeclassrealization}
Let $X$ be a prime Fano threefold of index 1 and degree 14. Let $\sigma'$ be a Serre-invariant stability condition on $\mathcal{A}_X$. Then the projection functor $\mathrm{pr}: \mathrm{D^b}(X)\rightarrow\mathcal{A}_X$ induces following isomorphisms between classical moduli spaces and Bridgeland moduli spaces of stable objects: 

\begin{enumerate}
    \item \emph{(Proposition \ref{conic-iso})} $\mathcal{C}(X)\xrightarrow{\cong} \mathcal{M}_{\sigma'}(\mathcal{A}_X, s)$.

    \item \emph{(Proposition \ref{216-iso})} $M_X(2,1,6)\xrightarrow{\cong} \mathcal{M}_{\sigma'}(\mathcal{A}_{X_{14}}, 3s-t)$.
    
    \item \emph{(Proposition \ref{2-16-iso})} $M_X(2,-1,6)\xrightarrow{\cong} \mathcal{M}_{\sigma'}(\mathcal{A}_{X_{14}}, 2s-t)$.
\end{enumerate}

\end{theorem}

\begin{theorem}\label{negative one class on cubic}
Let $Y$ be a cubic threefold. Let $\sigma$ be a Serre-invariant stability condition on $\mathcal{K}u(Y)$.Then the projection functor $\mathrm{pr}: \mathrm{D^b}(X)\rightarrow\mathcal{A}_X$ induces following isomorphisms between moduli spaces of stable sheaves and Bridgeland moduli spaces of stable objects: 

\begin{enumerate}

    \item \emph{(Proposition \ref{212-iso})} $M^b_Y(2,1,2)\xrightarrow{\cong} \mathcal{M}_{\sigma}(\mathcal{K}u(Y), v-w)$.
    
    \item \emph{(Proposition \ref{2-12-iso})} $M^b_Y(2,-1,2)\xrightarrow{\cong} \mathcal{M}_{\sigma}(\mathcal{K}u(Y), 2v-w)$.
\end{enumerate}

\end{theorem}

As a result, we recover two classical results:

\begin{corollary} \emph{(\cite[Theorem 7.2]{IM05})}
   Let $X$ be a prime Fano threefold of index 1 and  degree 14. Then we have following isomorphism:
   \[M_X(2,1,6)\cong \mathcal{C}(X)\]
\end{corollary}

\begin{corollary} \emph{(\cite[Theorem 1]{212cubic})}
Let $Y$ be a cubic threefold. Then we have following isomorphism:
\[M^b_Y(2,1,2)\cong \Sigma(Y)\]
\end{corollary}

Next we focus on the moduli spaces of instanton sheaves on index 2 prime Fano threefolds $Y:=Y_d$. Originally, instanton bundles appeared in \cite{adhm} as a way to describe Yang-Mills instantons on a 4-sphere $S^4$, which play an important role in Yang-Mills gauge theory. They provide extremely interesting links between physics and algebraic geometry. The notion of mathematical instanton bundle was first introduced on $\mathbb{P}^3$ and generalized to Fano threefolds in \cite{fae11} and \cite{kuz12}. On $Y=Y_3$ and $Y_4$, the instantonic condition is automatically satisfied for stable bundles of rank 2 with $c_1=0,c_2=2,c_3=0$ (for precise definitions, see Section \ref{sec-6}), as proved in \cite{D} and \cite{qin2019moduli}. However, this is not known on $Y_5$ and is conjectured in \cite[Conjecture 3.7]{qinv5}.


We denote the moduli space of S-equivalence classes of instanton sheaves by  $M^{inst}_Y$ and study the relations between the moduli space and the Bridgeland moduli spaces of S-equivalence classes of semistable objects of $(-4)$-classes in $\mathcal{K}u(Y_d)$. More precisely, we have:

\begin{theorem} \emph{(Theorem \ref{ins-iso}, Theorem \ref{ins-embd})} \label{1.6}
Let $Y:=Y_d$ be a prime Fano threefold of index 2 and degree $d$. When $d\neq 1,2$, the projection functor $\mathrm{pr}$ induces an isomorphism between moduli spaces of S-equivalence classes of instanton sheaves and Bridgeland moduli spaces:
\[M^{inst}_Y \xrightarrow{\cong}  \mathcal{M}^{ss}_{\sigma}(\mathcal{K}u(Y), 2-2L)\]
for every Serre-invariant stability condition $\sigma$ on $\mathcal{K}u(Y_d)$.

When $d=1,2$, the projection functor $\mathrm{pr}$ induces an open immersion:
\[M^{inst}_Y \hookrightarrow  \mathcal{M}^{ss}_{\sigma}(\mathcal{K}u(Y), 2-2L)\]
for every stability condition $\sigma\in \mathcal{K}$, where $\mathcal{K}$ is a certain family of stability conditions. 
\end{theorem}

In fact, one is able to show that there's an open immersion $M^{inst}_{Y_2} \hookrightarrow  \mathcal{M}^{ss}_{\sigma}(\mathcal{K}u(Y_2), 2-2L)$ for every Serre-invariant stability condition $\sigma$ on $\mathcal{K}u(Y_2)$.

On index one side, there're also some interesting moduli spaces of sheaves like moduli of instanton sheaves. As shown in \cite{BCFacm}, when $X$ is a non-hyperelliptic prime Fano threefold, there're some arithmetically Cohen-Macaulay (ACM) bundles in the moduli space $M_X^{ss}(2,0,4)$, and the sheaves in $M^{ss}_X(2,0,4)$ are classified when genus $g\geq 7$. Using the results in Section 7, we can show that:

\begin{theorem} \emph{(Theorem \ref{acm-iso})}
Let $X$ be a prime Fano threefold of index 1 and degree $14$, $18$ or $22$. Then the projection functor $\mathrm{pr}$ induces an isomorphism between moduli space of S-equivalence classes of  semistable sheaves and Bridgeland moduli space:
\[M^{ss}_X(2,0,4)\xrightarrow{\cong} \mathcal{M}^{ss}_{\sigma'}(\mathcal{A}_X, 2-4L)\]
such that the restriction gives
\[M_X(2,0,4)\xrightarrow{\cong} \mathcal{M}_{\sigma'}(\mathcal{A}_X, 2-4L)\]
for every Serre-invariant stability condition $\sigma'$ on $\mathcal{A}_X$.
\end{theorem}

As a corollary, we have following isomorphisms of moduli spaces of sheaves:

\begin{corollary}
Let $(Y,X)\in Z_d\subset \mathcal{MF}^2_{d}\times \mathcal{MF}^1_{4d+2}$ for $d=3,4,5$. Then we have following isomorphism of moduli spaces
\[M^{inst}_Y\cong M^{ss}_X(2,0,4)\]
\end{corollary}

\subsection{Related work}
The first example of stability conditions constructed in the Kuznetsov component of a cubic threefold is given in \cite{BMMS}. After a direct construction of stability conditions in the Kuznetsov components of a series Fano threefolds in \cite{BLMS}, the various Bridgeland moduli space of stable objects in the Kuznetsov components of index two and three Fano threefolds are studied in \cite{PY20}, \cite{APR}, \cite{petkovic2020note} \cite{bayer2020desingularization}, \cite{bolognese2021fullness}. In \cite{PY20}, the Fano surface of lines on $Y_d, d\neq 1$ is identified with moduli space of stable objects with $(-1)$-class $1-L$ in $\mathcal{K}u(Y_d)$, where  $d=1$ case is treated in \cite{petkovic2020note}. In \cite{zhang2020bridgeland} and \cite{JLZ}, the authors studied Bridgeland moduli spaces for Gushel-Mukai threefolds and used them to study several conjectures on them. In particular, the Bridgeland moduli spaces of stable objects with $(-1)$-classes in $\mathcal{A}_{X_{10}}$ are realized as $\mathcal{C}_m(X)$ and $M_X(2,1,5)$, while in our note, in Theorem~\ref{negativeclassrealization}, we  realize Bridgeland moduli space of stable objects of $(-1)$-classes in $\mathcal{A}_{X_{14}}$ as $\mathcal{C}(X)$, $M_X(2,1,6)$ and $M_X(2,1,-6)$. 

In our note, we identify moduli space of instanton sheaves with minimal charge on $Y_d, d=3,4,5$ with moduli space of (semi)stable objects of $(-4)$-class $2(1-L)$ 
in $\mathcal{K}u(Y_d)$. On cubic threefolds, these Bridgeland moduli spaces were studied in \cite{lahoz2015acm} via derived category of coherent sheaves of $\mathbb{P}^2$ with the action of a sheaf of Clifford algebras, and on $Y_4,Y_5$ these were studied in \cite{qin2019moduli} and \cite{qinv5} via classical stability of sheaves on curves and representations of quivers. 

When the first draft of our paper was finished, we learned that Xuqiang Qin independently proved similar results in 
\cite{qin2021bridgeland} to Theorem \ref{ins-iso}. He studied 
$M^{inst}_{Y_d}$ for $d=3,4,5$ and showed the isomorphism between $M^{inst}_{Y_d}$ and $\mathcal{M}^{ss}_{\sigma}(\mathcal{K}u(Y_d), 2-2L)$ 
for a family of stability conditions $\mathcal{K}$ constructed in \cite{PY20}, while we show this for every Serre invariant stability condition defined on $\mathcal{K}u(Y_d)$. Thus, we use slightly different techniques in proving Theorem \ref{ins-iso}.  In \cite{qin2021bridgeland}, the author applies wall-crossing techniques developed in \cite{PY20} and \cite{bayer2020desingularization}, while we have to not only use wall-crossing techniques but also use weak Mukai lemma to prove stability for all Serre-invariant stability conditions. He also showed the stability of instanton bundles of minimal charge, while we prove this for all instanton sheaves in Proposition \ref{sheaf-stable}. Moreover, in \cite{qin2021bridgeland}, the author makes use of the notion of 2-Gieseker stability while we do not but apply more elementary techniques. In addition we prove Lemma \ref{same-slope}, Lemma \ref{non-negative} and Proposition \ref{all-serre-invariant} to make the theorem work for every Serre-invariant stability condition.

\subsection{Further questions}
 It is an interesting question to classify instanton sheaves on a general quartic double solid $Y_2$ and prove similar result as Theorem~\ref{ins-iso}. On the other hand, in \cite{zhang2020bridgeland}, it is shown that $\mathcal{K}u(Y_2)$ is not equivalent to $\mathcal{A}_{X_{10}}$, so it is interesting to study the relation between Bridgeland moduli spaces $\mathcal{M}_{\sigma}(\mathcal{K}u(Y_2),2-2L)$ and $\mathcal{M}_{\sigma'}(\mathcal{A}_{X_{10}},2-4L)$ for Serre invariant stability conditions.


\subsection{Structure of the paper}

In Section \ref{sec-2} and Section \ref{sec-3}, we recall some basic definitions and properties of Kuznetsov component  and (weak)stability conditions. Fix Serre-invariant stability conditions $\sigma$ and $\sigma'$ on $\mathcal{K}u(Y)$ and $\mathcal{A}_X$ respectively.  In Section \ref{sec-4}, we first show that all moduli spaces of stable objects of $(-1)$-classes in $\mathcal{K}u(Y)$ and $\mathcal{A}_X$ are isomorphic. Then we show that the projection functor $\mathrm{pr}$ induces isomorphisms $\mathcal{C}(X)\cong\mathcal{M}_{\sigma'}(\mathcal{A}_X, s)$ for every $X:=X_{4d+2}$ and $d=3,4,5$. In Section \ref{sec-5} we focus on $Y_3$ and $X_{14}$. Using similar arguments in Section \ref{sec-4}, we realize three moduli spaces of stable objects of $(-1)$-classes as classical moduli space. In Section \ref{sec-6} we first recall some definitions and properties of instanton sheaves, and show that the instantonic condition is closely related to the walls of Bridgeland stability conditions.
In Section \ref{sec-7} we use the results in Section \ref{sec-6} to show that for $Y:=Y_d$ the projection functor $\mathrm{pr}$ gives an isomorphism $M^{inst}_Y\to  \mathcal{M}^{ss}_{\sigma}(\mathcal{K}u(Y), 2v)$ when $d\neq 1,2$ and open immersions when $d=1,2$. In Section \ref{sec-8} we show that for every pair $(Y_d, X_{4d+2})\in Z_d\subset \mathcal{MF}^2_{d}\times \mathcal{MF}^1_{4d+2}$ and $d=3,4,5$, the equivalence $\mathcal{K}u(Y_d)\cong \mathcal{A}_{X_{4d+2}}$ gives an isomorphism $\mathcal{M}_{\sigma}(\mathcal{K}u(Y), 2v)\cong \mathcal{M}_{\sigma'}(\mathcal{A}_X, 2s)$. As a corollary, when $d=3,5$ we show that all moduli spaces of semistable objects of $(-4)$-classes in $\mathcal{K}u(Y)$ and $\mathcal{A}_X$ are isomorphic. As an application, in Section \ref{sec-9} we show that for $X:=X_{4d+2}$ and $d=3,4,5$ the projection functor gives an isomorphism $M^{ss}_X(2,0,4)\cong \mathcal{M}^{ss}_{\sigma'}(\mathcal{A}_X, 2s)$.

\subsection{Acknowledgements}
This note originated from a seminar talk given by the second author in Tianyuan Mathematical Center in Southwest China(TMCSC), Sichuan University. It is part of the Undergraduate Mathematical Research Project of the first author mentored by the second author. We would like to thank professor Xiaojun Chen for the invitation and TMCSC for their hospitality. We thank Daniele Faenzi, Li Lai, Girivaru Ravindra, Junyan Xu, and Song Yang for answering us several questions. We also thank Arend Bayer, Augustinas Jacovskis for useful conversations on several related topics. The first author would like to thank Jiahui Gao, Songtao Ma, Rui Xiong, and Jiajin Zhang for useful discussion. The second author thanks Tingyu Sun for support. 
In addition, we would like to thank Laura Pertusi and Xuqiang Qin for informing us the very recent preprint \cite{qin2021bridgeland} and giving useful comments on the earlier draft of our paper.

\section{Preliminaries} \label{sec-2}

\subsection{Prime Fano threefolds}

A complex smooth projective variety with ample anticanonical bundle is called Fano. A Fano variety is called prime if it has Picard number 1. For a prime Fano variety $X$, we can choose a unique ample divisor such that $\mathrm{Pic}(X)\cong \mathbb{Z}\cdot H$, which is called the fundamental divisor of $X$. The index of a prime Fano variety is the least integer $i$ such that $-K_X=i\cdot H$. The degree of a prime Fano variety is $d:=H^3$. The genus $g$ is defined as $2g-2=d$.

\subsection{Derived Category of Fano Threefolds}

Let $Y_d$ be a prime Fano threefold of index 2 and degree d. We consider the semiorthogonal decomposition of $Y_d$ given by 
\[\mathrm{D^b}(Y_d)=\langle \mathcal{K}u(Y_d), \mathcal{O}_{Y_d}, \mathcal{O}_{Y_d}(1)\rangle\]
The numerical Grothendieck group $\mathcal{N}(\mathcal{K}u(Y_d))$ of $\mathcal{K}u(Y_d)$ is a rank 2 lattice spanned by 
\[v:=[I_L]=1-\frac{1}{d}H^2, w:=H-\frac{1}{2}H^2+(\frac{1}{6}-\frac{1}{d})H^3\]
where $I_L$ is an ideal sheaf of a line $L$ on $Y_d$. The Euler form is given by 
\begin{equation}
	 \left[               
	\begin{array}{cc}   
	-1 & -1 \\  
	1-d & -d \\
	\end{array}
	\right] 
\end{equation}
In the case of index one, the semiorthogonal decompositions of $X_{2g-2}$ of even genus $6\leq g\leq 12$ are given in \cite{kuz09} and \cite{kuz06}.
\[\mathrm{D^b}(X_{2g-2})=\langle\mathcal{K}u(X_{2g-2}), \mathcal{E}, \mathcal{O}_{X_{2g-2}}\rangle\]
where $\mathcal{E}$ is a rank $2$ stable vector bundle with ch$(\mathcal{E})=2-H+\frac{g-4}{2}L+\frac{10-g}{12}P$. However we will use another semiorthogonal decomposition in \cite[Section B.2]{KPS}, which is given by
\[\mathrm{D^b}(X_{2g-2})=\langle\mathcal{A}_{X_{2g-2}}, \mathcal{O}_{X_{2g-2}}, \mathcal{E}^{\vee}\rangle\]

The  numerical Grothendieck group $\mathcal{N}(\mathcal{A}_{X_{2g-2}})$ of $\mathcal{A}_{X_{2g-2}}$ is a rank 2 lattice spanned by
\[s:=[I_C]=1-\frac{1}{g-1}H^2, t:=H-(\frac{g}{2}+1)L-\frac{16-g}{12}P\]
where $I_C$ is an ideal sheaf of a conic $C$ on $X_{2g-2}$. The Euler form is given by
\begin{equation}
\left[               
\begin{array}{cc}   
-1 & -2 \\  
1-\frac{g}{2} & 1-g \\
\end{array}
\right] 
\end{equation}


	
		

As noted in \cite[Proposition 3.9]{kuz09}, the Chern character map identify the numerical Grothendieck group with the lattice generated by Chern characters of some sheaves, then we will use notations of numerical class and Chern character alternatively.

Up to sign, there're three $(-1)$-classes $s,3s-t$ and $2s-t$ in $\mathcal{A}_{X_{14}}$, and three $(-1)$-classes $v,2v-w$ and $v-w$ in $\mathcal{K}u(Y_3)$.
And there're three $(-4)$-classes, which are the double of three $(-1)$-classes.

It's also easy to see that there're infinitely many $(0)$-class, $(-1)$-class and $(-4)$-class, and no $(-2)$-class and $(-3)$-class in $\mathcal{K}u(Y_{4})\cong \mathcal{A}_{X_{18}}$; And there're infinitely many $(-1)$-class and $(-4)$-class, and no $(0)$-class, $(-2)$-class and $(-3)$-class in $\mathcal{K}u(Y_{5})\cong \mathcal{A}_{X_{22}}$.


We also need some results from \cite{mu92}, \cite{KPS} and \cite{BCFacm}:

\begin{lemma} \label{215}
	Let $\mathrm{D^b}(X_{2g-2})=\langle\mathcal{A}_{X_{2g-2}}, \mathcal{O}_X, \mathcal{E}^{\vee}\rangle$ as above, then
	
	\begin{enumerate}
	    \item For each $g=8,10,12$ there's a closed immersion $X_{2g-2}\hookrightarrow \mathrm{Gr}(2,\frac{g}{2}+2)$ such that $\mathcal{E}$ is the pullback of the tautological bundle on $\mathrm{Gr}(2,\frac{g}{2}+2)$.
	
		\item $\mathcal{E}^{\vee}$ is the unique rank $2$ stable vector bundle with $c_1=H, c_2=\frac{g+2}{2}L, c_3=0$, called Mukai bundle; Moreover, $\mathcal{E}^{\vee}$ is globally generated and ACM.
		
		\item 
		 $\mathrm{hom}(\mathcal{E}^{\vee}, \mathcal{E}^{\vee})=1, \mathrm{ext}^i(\mathcal{E}^{\vee}, \mathcal{E}^{\vee})=0, ~\forall i\geq 1$ and $h^i(X, \mathcal{E}^{\vee}(-1))=0$ for all $i$.
	\end{enumerate}
\end{lemma}

\subsection{Moduli spaces of sheaves}

In this subsection, we recall some definitions about (semi)stable sheaves. We refer to book \cite{HL} for a more detailed account of notions.

Let $(X, H)$ be a smooth polarized $n$-dimensional projective variety. Recall that a torsion-free sheaf $F$ is \emph{Gieseker-semistable} if for any coherent subsheaf $E$ with $0<\mathrm{rk}(E)<\mathrm{rk}(F)$, one has $\frac{p(E,t)}{\mathrm{rk}(E)}\leq \frac{p(F,t)}{\mathrm{rk}(F)}$ for $t\gg 0$. The sheaf $F$ is called \emph{Gieseker-stable} if the inequality above is always strict.

The \emph{slope} of a sheaf $F$ with positive rank is defined as $\mu(F):=\frac{c_1(F)H^{n-1}}{\mathrm{rk}(F)H^n}$. We recall that a torsion-free coherent sheaf $F$ is \emph{$\mu$-semistable} if for any coherent subsheaf $E$ with $0<\mathrm{rk}(E)<\mathrm{rk}(F)$, on has $\mu(E)\leq \mu(F)$. The sheaf $F$ is called \emph{$\mu$-stable} if the above inequality is always strict. The two notions are related as following:
\[\mu-\text{stable} \Rightarrow \text{Gieseker stable} \Rightarrow \text{Gieseker semistable} \Rightarrow \mu-\text{semistable}\]
We denote the moduli space of $S$-equivalence classes of rank $r$ torsion-free Gieseker-semistable sheaves with Chern classes $c_1, c_2$ and $c_3=0$ by $M^{ss}_X(r,c_1,c_2)$. And we denote the moduli space of $S$-equivalence classes of rank $r$ torsion-free Gieseker-stable sheaves with Chern classes $c_1, c_2$ and $c_3=0$ by $M_X(r,c_1,c_2)$.

\section{Review on Bridgeland stability conditions} \label{sec-3}

In this section, we recall the definition and some properties of (weak)stability conditions on triangulated category.

\subsection{(Weak)stability conditions}

Let $\mathcal{T}$ be a triangulated category.

\begin{definition}
	A heart of a bounded t-structure on $\mathcal{T}$ is a full subcategory $\mathcal{A}\subset \mathcal{T}$ such that 
	
	\begin{enumerate}
		\item for $A,B\in \mathcal{A}$ and $n<0$ we have $\mathrm{Hom}(A,B[n])=0$, and
		
		\item for every object $F\in \mathcal{T}$ there exists a sequence of morphisms
		\[0=F_0\xrightarrow{\phi_1} F_1\to ... \xrightarrow{\phi_m} F_m=F\]
		such that cone$(\phi_i)$ is of form $A_i[k_i]$ for some sequence $k_1>k_2>...>k_m$ of integers and objects $A_i\in \mathcal{A}$. We denote $A_i$ by $\mathcal{H}^{-k_i}(F)$.
	\end{enumerate}
\end{definition}

\begin{definition}
	Let $\mathcal{A}$ be an abelian category. A \emph{weak stability function} on $\mathcal{A}$ is a group homomorphism
	$Z: K(\mathcal{A})\to  \mathbb{C}$ such that for every non-zero object $A\in \mathcal{A}$, we have
	$\mathrm{Im}Z(A)\geq 0$, and $\mathrm{Im}Z(A)=0\Rightarrow \mathrm{Re}Z(A)\leq 0$.
	We say $Z$ is a \emph{stability function} if $\mathrm{Im}Z(A)\geq 0, ~\text{and} ~\mathrm{Im}Z(A)=0\Rightarrow \mathrm{Re}Z(A)< 0$.
\end{definition}

Fix a finite rank lattice $\Lambda$ and a surjective homomorphism $v: K(\mathcal{A})\twoheadrightarrow \Lambda$.

\begin{definition}
	A \emph{weak stability condition on $\mathcal{T}$} with respect to the lattice $\Lambda$ is a pair $\sigma=(\mathcal{A}, Z)$, where $Z: \Lambda \to \mathbb{C}$ is a group homomorphism and $\mathcal{A}$ is the heart of a bounded t-structure, satisfies the following conditions:
	
	\begin{enumerate}
		\item The composition 
		$K(\mathcal{A})= K(\mathcal{T})\stackrel{v}{\to} \Lambda \stackrel{Z}{\to}\mathbb{C}$
		 is a weak stability function on $\mathcal{A}$; we denote $Z(A):=Z(v(A))$. We define a \emph{slope} for any object $A\in \mathcal{A}$ as $$\mu_{\sigma}(A):=-\frac{\mathrm{Re}Z(A)}{\mathrm{Im}Z(A)}$$ for $\mathrm{Im}Z(A)>0$ and $\mu_{\sigma}(A):=+\infty$ otherwise.
		 
		 \item A non-zero object $A\in \mathcal{A}$ is called \emph{$\sigma$-semistable} (resp. \emph{$\sigma$-stable}) if for every non-zero proper subobject $B\subset A$, we have $\mu_{\sigma}(B)\leq \mu_{\sigma}(A)$ (resp. $\mu_{\sigma}(B)< \mu_{\sigma}(A)$). An object $F\in \mathcal{T}$ is called  $\sigma$-(semi)stable if $F[k]\in \mathcal{A}$ is $\sigma$-(semi)stable for some $k\in \mathbb{Z}$.
		 
		 \item Every object $A\in \mathcal{A}$ has a Harder-Narasimhan filtration in $\sigma$-semistble objects.
		 
		 \item There is a quadratic form $Q$ on $\Lambda\otimes \mathbb{R}$ such that $Q|_{\mathrm{Ker}(Z)}$ is negative definite, and $Q(A)\geq 0$ for all $\sigma$-semistable $A\in \mathcal{A}$.
	\end{enumerate}
\end{definition}

A weak stability condition $\sigma=(\mathcal{A}, Z)$ on $\mathcal{T}$ with respect to the lattice $\Lambda$ is called a \emph{Bridgeland stability condition} if $Z\circ v$ is a stability function on $\mathcal{A}$.

\begin{definition}
	The \emph{phase} of a $\sigma$-semistable object $A\in \mathcal{A}$ is 
	\[\phi(A):=\frac{1}{\pi} \mathrm{arg}(Z(A))\in (0, 1]\]
	and for $A[n]$, we set 
	$\phi(A[n]):=\phi(A)+n.$
	
	A \emph{slicing} $\mathcal{P}_{\sigma}$ of $\mathcal{T}$ is a collection of full additive subcategories $\mathcal{P}_{\sigma}(\phi)\subset \mathcal{T}$ for $\phi \in \mathbb{R}$, such that the subcategory $\mathcal{P}_{\sigma}(\phi)$ is given by the zero object and all $\sigma$-semistable objects with phase $\phi$.
	
\end{definition}

We both use the notation $\sigma=(\mathcal{P}_{\sigma}, Z)$ and $\sigma=(\mathcal{A}_{\sigma}, Z)$ for a (weak)stability condition with heart $\mathcal{A}_{\sigma}:=\mathcal{P}_{\sigma}((0,1])$, where $\mathcal{P}_{\sigma}$ is a slicing and $\mathcal{P}_{\sigma}((0,1])$ is the extension closure of subcategories $\{\mathcal{P}_{\sigma}(\phi) ~|~ \phi \in (0,1]\}$.

We denote by Stab$(\mathcal{T})$ the set of stability conditions on $\mathcal{T}$. The universal covering space $\tilde{\mathrm{GL}^+_2(\mathbb{R})}$ of $\mathrm{GL}^+_2(\mathbb{R})$ has a right action on Stab$(\mathcal{T})$, defined in \cite[Lemma 8.2]{T03}.

Starting from a weak stability condition $\sigma=(\mathcal{A}, Z)$ on $\mathcal{T}$, we can produce a new heart of a bounded t-structure, by \emph{tilting} $\mathcal{A}$. Let $\mu\in \mathbb{R}$, we define following subcategoies of $\mathcal{A}$:
\[\mathcal{T}^{\mu}_{\sigma}:=\{A\in \mathcal{A}: ~ \text{all HN factors}~ B ~\text{of}~A~ \text{have slope}~\mu_{\sigma}(B)>\mu\}\]
\[\mathcal{F}^{\mu}_{\sigma}:=\{A\in \mathcal{A}: ~ \text{all HN factors}~ B ~\text{of}~A~ \text{have slope}~\mu_{\sigma}(B)\leq \mu\}\]
Thus by \cite{HRS}, the category
$\mathcal{A}^{\mu}_{\sigma}:=\langle \mathcal{T}^{\mu}_{\sigma}, \mathcal{F}^{\mu}_{\sigma}[1]\rangle$
is the heart of a bounded t-structure on $\mathcal{T}$.

We say that the heart $\mathcal{A}^{\mu}_{\sigma}$ is obtained by tilting $\mathcal{A}$ with respect to the weak stability condition $\sigma$ at the slope $\mu$.

\subsection{Weak stability condition on $\mathrm{D^b}(X)$}

Let $X$ be a smooth projective variety of dimension $n$ and $H$ be an ample divisor on $X$. Following \cite[Section 2]{BLMS}, we review the construction of weak stability conditions on $\mathrm{D^b}(X)$.

For any $j\in \{0,1,2,...,n\}$, we consider the lattice $\Lambda^j_H\cong \mathbb{Z}^{j+1}$ generated by 
\[(H^n\text{ch}_0, H^{n-1}\text{ch}_1, ..., H^{n-j}\text{ch}_j)\in \mathbb{Q}^{j+1}\]
with surjective map $v^j_H: K(X)\to \Lambda^j_H$ induced by Chern character. 
The pair 
$\sigma_H:=(\mathrm{Coh}(X), Z_H)$
where $Z_H: \Lambda^1_H\to \mathbb{C}$ is given by 
$Z_H(F):=-H^{n-1}\text{ch}_1(F)+iH^n \text{ch}_0(F)$
defines a weak stability condition on $\mathrm{D^b}(X)$ with respect to the lattice $\Lambda^1_H$ (see \cite[Example 2.8]{BLMS}). 
Moreover, any $\mu_H$-semistable sheaf $F$ staifies the following Bogomolov-Gieseker inequality:
$$\Delta_H(F):=(H^{n-1}\text{ch}_1(F))^2-2H^n\text{ch}_0(F)\cdot H^{n-2}\text{ch}_2(F)\geq 0.$$

Given a parameter $\beta \in \mathbb{R}$, we denote by
$\mathrm{Coh}^{\beta}(X)$
the heart of bounded t-structure obtained by tilting the weak stability condition $\sigma_H$ at the slope $\mu_H=\beta$. For $F\in \mathrm{D^b}(X)$, we set 
$\text{ch}^{\beta}(F):=e^{-\beta}\text{ch}(F)$.
Explicitly, for each non-negative integer $k$ we have 
$$\mathrm{ch}^{\beta}_{k}(F)=\sum_{i=0}^k \frac{(-\beta)^iH^i}{i!}\mathrm{ch}_{k-i}(F).$$

\begin{proposition} \emph{(\cite[Proposition 2.12]{BLMS})}. \label{blms} For any  $(\alpha, \beta)\in \mathbb{R}_{>0}\times \mathbb{R}$, the pair 
$\sigma_{\alpha, \beta}=(\mathrm{Coh}^{\beta}(X), Z_{\alpha, \beta})$
with
\[Z_{\alpha, \beta}(F):=\frac{1}{2}\alpha^2H^n\mathrm{ch}^{\beta}_0(F)-H^{n-2}\mathrm{ch}^{\beta}_2(F)+iH^{n-1}\mathrm{ch}^{\beta}_1(F)\]
defines a weak stability condition on $\mathrm{D^b}(X)$ with respect to $\Lambda^2_H$. The quadractic form $Q$ can be given by the discriminant $\Delta_H$. Moreover, these weak stability conditions vary continuously as $(\alpha, \beta)\in \mathbb{R}_{>0}\times \mathbb{R}$ varies.
\end{proposition}

We can visualize the weak stability conditions $\sigma_{\alpha, \beta}$ in the upper half plane $\mathbb{R}_{>0}\times \mathbb{R}$.

\begin{definition}
	Let $v$ be a vector in $\Lambda^2_H.$
	
	\begin{enumerate}
		\item A \emph{numerical wall} for $v$ is the set of pairs $(\alpha, \beta)\in \mathbb{R}_{>0}\times \mathbb{R}$ such that there is a vector $w\in \Lambda^2_H$ verifying the numerical relation $\mu_{\alpha, \beta}(v)=\mu_{\alpha, \beta}(w)$.
		
		\item A \emph{wall} for $E\in \mathrm{Coh}^{\beta}(X)$ is a numerical wall for $v:=\text{ch}_{\leq 2}(E)$, where $\text{ch}_{\leq 2}(E):=(\text{ch}_0(E), \text{ch}_1(E), \text{ch}_2(E))$, such that for every $(\alpha, \beta)$ on the wall there's an exact sequence of semistable objects $0\to F\to E\to G\to 0$ in $\mathrm{Coh}^{\beta}(X)$ such that $\mu_{\alpha, \beta}(F)=\mu_{\alpha, \beta}(E)=\mu_{\alpha, \beta}(G)$ gives rise to the numerical wall.
		
		\item A \emph{chamber} is a connected component in the complement of the union of walls in the upper half-plane.
	\end{enumerate}
\end{definition}

An important property is that the weak stability conditions $\sigma_{\alpha, \beta}$ satisfy well-behaved wall-crossing: walls respect to a class $v\in \Lambda^2_H$ are locally finite. By \cite[Proposition B.5]{BMS}, if $v=\text{ch}_{\leq 2}(E)$ with $E\in \mathrm{Coh}^{\beta}(X)$, then the stability of $E$ remains unchanged as $(\alpha, \beta)$ varies in chamber.

Finally, we recall following variant of the weak stability conditions of Proposition \ref{blms}, which will be used frequently. Fix $\mu\in \mathbb{R}$ and let $u$ be the unit vector in upper half plane with $\mu=-\frac{\mathrm{Re}u}{\mathrm{Im}u}$. We denote by
$\mathrm{Coh}^{\mu}_{\alpha, \beta}(X)$
the heart obtained by tilting the weak stability condition $\sigma_{\alpha, \beta}=(\mathrm{Coh}^{\beta}(X), Z_{\alpha, \beta})$ at the slope $\mu_{\alpha, \beta}=\mu$.

\begin{proposition} \emph{\cite[Proposition 2.15]{BLMS}}. The pair $\sigma^{\mu}_{\alpha, \beta}:=(\mathrm{Coh}^{\mu}_{\alpha, \beta}(X), Z^{\mu}_{\alpha, \beta})$, where
\[Z^{\mu}_{\alpha, \beta}:=\frac{1}{u}Z_{\alpha, \beta}\]
is a weak stability condition on $\mathrm{D^b}(X)$.
\end{proposition}

\subsection{Serre-invariant stability conditions on Kuznetsov component}

In this subsection, we recall some results in \cite{PY20}, which we will use frequently in the next sections.

\begin{definition}
	A stability $\sigma$ on $\mathcal{T}$ is called \emph{Serre-invariant} if $S_{\mathcal{T}}(\sigma)=\sigma\cdot\tilde{g}$ for some $\tilde{g}\in \tilde{\mathrm{GL}^+_2(\mathbb{R})}$.
\end{definition}

Recall that $\mathcal{K}u(Y_4)\cong \mathcal{A}_{X_{18}}\cong \mathrm{D^b}(C_2)$ and $\mathcal{K}u(Y_5)\cong \mathcal{A}_{X_{22}}\cong \mathrm{D^b}(Q_3)$. Here $C_2$ is a smooth curve of genus 2 and $Q_3$ is the 3-Kronecker quiver (see \cite{KPS}). Stability conditions on these two categories are studied in \cite{m07} and \cite{kronecker}.

Following \cite{PY20}, we define
\[V:=\{(\alpha, \beta)\in \mathbb{R}_{>0}\times \mathbb{R} ~|~ -\frac{1}{2}\leq \beta <0, \alpha<-\beta, ~ \text{or} ~ -1<\beta <-\frac{1}{2}, \alpha\leq 1+\beta\}\]
and let $\mathcal{K}$ be the orbit of $V$ by $\tilde{\mathrm{GL}^+_2(\mathbb{R})}$. We set $Z(\alpha, \beta):=Z^0_{\alpha, \beta}|_{\mathcal{K}u(Y)}$, and $\mathcal{A}(\alpha, \beta):=\mathrm{Coh}^0_{\alpha, \beta}(Y)\cap \mathcal{K}u(Y)$. We define a lattice
$\Lambda^2_{H, \mathcal{K}u(Y)}:=Im(K(\mathcal{K}u(Y))\to K(Y)\to \Lambda^2_H)\cong \mathbb{Z}^2$.

\begin{theorem} \emph{(\cite{BLMS}, \cite[Theorem 3.3, Proposition 3.6, Corollary 5.5]{PY20})}  Let $Y$ be a prime Fano threefold of index $2$.
	
	\begin{enumerate}
		\item  The pair 
		$\sigma(\alpha, \beta):=(\mathcal{A}(\alpha, \beta), Z(\alpha, \beta))$
		is a Bridgeland stability condition on $\mathcal{K}u(Y)$ with respect to $\Lambda^2_{H, \mathcal{K}u(Y)}\cong \mathbb{Z}^2$ for every $(\alpha, \beta)\in V$.
		
		\item Fix $0<\alpha_0<\frac{1}{2}$. For every $(\alpha, \beta)\in V$, there's a $\tilde{g}\in \tilde{\mathrm{GL}^+_2(\mathbb{R})}$ such that $\sigma(\alpha, \beta)=\sigma(\alpha_0, -\frac{1}{2})\cdot \tilde{g}$. Hence $\mathcal{K}=V\cdot \tilde{\mathrm{GL}^+_2(\mathbb{R})}=\sigma(\alpha_0, -\frac{1}{2})\cdot\tilde{\mathrm{GL}^+_2(\mathbb{R})}$.
		
		\item Every stability condition in $\mathcal{K}$ is Serre-invariant.
	\end{enumerate}
	
\end{theorem}

Thus via the equivalences $\mathcal{K}u(Y_d)\cong \mathcal{A}_{X_{4d+2}}$ for $d=3,4,5$, since this equivalence commutes with Serre functor, a Serre-invariant stability condition on $\mathcal{K}u(Y_d)$ induces a Serre-invariant stability condition on $\mathcal{A}_{X_{4d+2}}$.

As shown in \cite[Section 5.2]{PY20}, there're some useful numerical properties for Serre-invariant stability condition on $\mathcal{K}u(Y)$. Note that in the proof of following lemmas in \cite{PY20}, the only properties of $Y$ they used is that $S^3_{\mathcal{K}u(Y)}\cong [5]$, hence via the equivalence $\Phi$ in Theorem \cite{KPS}, these properties also hold for Serre-invariant stability conditions on $\mathcal{A}_X$.

For $X=X_{14}$, we have:

\begin{lemma} 
\label{py-phase}
	For every Serre-invariant stability condition $\sigma$ on $\mathcal{A}_{X_{14}}$, if $F\in \mathcal{A}_{X_{14}}$ is a $\sigma$-semistable object of phase $\phi(F)$, then the phase of $S_{\mathcal{A}_{X_{14}}}(F)$ satisfies $\phi(F)< \phi(S_{\mathcal{A}_{X_{14}}}(F))< \phi(F)+2$. In particular, we have $\mathrm{Ext}^2(F,F)=0$.
\end{lemma}

\begin{lemma}
\label{heart-dim}
	The heart of a Serre-invariant stability condition  $\sigma$ on $\mathcal{A}_{X_{14}}$ has homological dimension 2.
\end{lemma}

\begin{lemma} 
\label{py-ext-2}
	 For every Serre-invariant stability condition $\sigma=(\mathcal{A}, Z)$ on $\mathcal{A}_{X_{14}}$ and every non-zero object $F\in \mathcal{A}$, we have $\chi(F,F)\leq -1$ and  $\mathrm{ext}^1(F,F)\geq 2$.
\end{lemma}

\begin{lemma} 
\label{mk-lem} \emph{(Weak Mukai Lemma)}
	Let $\sigma$ be a Serre-invariant stability condition on $\mathcal{A}_{X_{2g-2}}$. Let $A\to F\to B$ be a triangle in the heart  $\mathcal{A}_{\sigma}$ such that $\mathrm{hom}(A,B)=\mathrm{Hom}(B, A[2])=0$. Then
	\[\mathrm{ext}^1(A,A)+\mathrm{ext}^1(B,B)\leq \mathrm{ext}^1(F,F)\]
\end{lemma}

When $g=8$, as in \cite{PY20}, $\mathrm{Hom}(B,A[2])=0$ if $\sigma$-semistable factors of $A$ have phases greater or equal than the phases of the $\sigma$-semistable factors of $B$.

When $g=10$ and $12$, since every heart of a stability condition on $\mathcal{A}_{2g-2}$ has homological dimension 1, hence automatically we have $\mathrm{Hom}(B, A[2])=0$.

The same argument in \cite[Lemma 5.13]{PY20} shows that:

\begin{lemma} 
\label{ku-stable}
	Let $X:=X_{14}$ or $X_{22}$. Let $\sigma$ be a Serre-invariant stability condition on $\mathcal{A}_X$. Then every $F\in \mathcal{A}_X$ with $\mathrm{ext}^1(F,F)=2$ is $\sigma$-stable. If $X=X_{18}$, then every $F\in \mathcal{A}_X$ with $\mathrm{ext}^1(F,F)=1$ or $\chi(F,F)=-1$ and $\mathrm{ext}^1(F,F)=2$ is $\sigma$-stable.
\end{lemma}


\begin{theorem} \emph{(\cite[Theorem 1.2]{PY20})} \label{sm-moduli} Let $X:=X_{14}, X_{18}$ or $X_{22}$. Every non-empty moduli space of $\sigma'$-stable objects in $\mathcal{A}_X$ with respect to a Serre-invariant stability condition $\sigma'$ is smooth.
\end{theorem}


       

\section{Isomorphisms between Bridgeland moduli spaces of $(-1)$-classes} \label{sec-4}

In this section, we fix $(Y, X)=(Y_d, X_{4d+2})\in Z_d\subset  \mathcal{MF}^2_{d}\times \mathcal{MF}^1_{4d+2}$ as in \cite{KPS} for $d=3,4,5$.
We denote the closed point in $\mathcal{M}_{\sigma'}$ that corresponds to $E\in \mathcal{A}_X$ by $[E]$.



\begin{proposition} \label{iso-moduli}
	The equivalence $\Phi: \mathcal{K}u(Y)\xrightarrow{\cong} \mathcal{A}_{X}$ induces an isomorphism between moduli spaces
	\[\mathcal{M}_{\sigma}(\mathcal{K}u(Y), 1-L)\xrightarrow{\cong} \mathcal{M}_{\sigma'}(\mathcal{A}_{X}, 1-2L)\]
	such that maps $[E]$ to $[\Phi(E)]$ on the level of closed points. Here $\sigma, \sigma'$ are Serre-invariant stability conditions on $\mathcal{K}u(Y), \mathcal{A}_X$ respectively. 
	In particular, $\mathcal{M}_{\sigma'}(\mathcal{A}_{X}, 1-2L)$ is irreducible and smooth of dimension 2.
\end{proposition}

\begin{proof}
	
	From \cite[Theorem 1.1]{PY20} we know that $\Sigma(Y)\cong \mathcal{M}_{\sigma}(\mathcal{K}u(Y), 1-L)$, hence there's a universal family on $\mathcal{M}_{\sigma}(\mathcal{K}u(Y), 1-L)$. Since $\Phi$ is of Fourier-Mukai type, the same argument in Lemma \ref{lem-4-1} shows that $\Phi$ induces a morphism such that maps $[E]$ to $[\Phi(E)]$. By Lemma \ref{ku-stable} we know that this morphism is a bijection between closed points.
	Now since $\Phi$ is an equivalence, the induced morphism is etale. Therefore, this morphism is bijective and etale, which is an isomorphism. The last statement also follows from \cite[Theorem 1.1]{PY20}.
	
\end{proof}

\begin{corollary}\label{all moduli spaces of (-1)-class isom}
Let $a$ be a $(-1)$-class in $\mathcal{N}(\mathcal{K}u(Y))$ and $a'$ be a $(-1)$-class in $\mathcal{N}(\mathcal{A}_{X})$. Then the Bridgeland moduli spaces $\mathcal{M}_{\sigma}(\mathcal{K}u(Y),a)\cong\mathcal{M}_{\sigma'}(\mathcal{A}_{X},a')$ for any Serre invariant stability condition $\sigma,\sigma'$ on $\mathcal{K}u(Y), \mathcal{A}_X$ respectively. 
\end{corollary}

\begin{proof}

\begin{enumerate}
    \item When $d=3$, all $(-1)$-classes in $\mathcal{A}_{X_{14}}$ up to signs are $[I_C]=1-2L, [S_{\mathcal{A}_X}(I_C)]$ and  $[S_{\mathcal{A}_X}^2(I_C)]$, and all $(-1)$-classes in $\mathcal{K}u(Y_3)$ up to signs are $[I_L]=1-L$, $[S_{\mathcal{K}u(Y_3)}(I_L)]$ and $[S^2_{\mathcal{K}u(Y_3)}(I_L)]$. Now $\sigma$ and $\sigma'$ are Serre-invariant stability conditions, thus Bridgeland moduli spaces of stable objects with these $(-1)$-classes are all isomorphic. 
    
    \item When $d=4$, there are infinitely many $(-1)$-classes in $\mathcal{N}(\mathcal{K}u(Y_4))$ and $\mathcal{N}(\mathcal{A}_{X_{18}})$ respectively. Note that $\mathcal{K}u(Y_4)\cong\mathcal{A}_{X_{18}}\cong D^b(C_2)$. All $(-1)$-classes in $\mathcal{K}u(Y_4)$ are permuted by rotation functor $R:E\mapsto\mathrm{L}_{\mathcal{O}_Y}(E\otimes\mathcal{O}_Y(H))$ and the functor preserves the stability conditions on $\mathcal{K}u(Y_4)$ by \cite[Proposition 5.7]{PY20}. Then by similar arguments in Proposition~\ref{iso-moduli}, all the moduli spaces are isomorphic. 
    
    \item $d=5$. Let $xv+yw\in\mathcal{N}(\mathcal{K}u(Y_5))$ be a $(-1)$-class, then the pair of integers $(x,y)$ are solutions of equation $x^2+5xy+5y^2=1$. Up to a linear transform, it is the Pell's equation $x^2-5y^2=1$.Up to sign, all the integer solutions $(x,y)$ are parametrized by $\eta^i,i\in\mathbb{Z}$, where $\eta=\frac{1}{2}(3+\sqrt{5})$ is the fundamental unit. On the other hand, the rotation functor $R:\mathcal{K}u(Y_5)\rightarrow\mathcal{K}u(Y_5)$ acts on $(-1)$-classes by $(x,y)\mapsto (-4x-5y,x+y)$. It is easy to check that the action of $R$ is one to one corresponds to the multiplication by $\eta$. This means that all $(-1)$-classes in $\mathcal{N}(\mathcal{K}u(Y_5))$ are given by action of rotation functor $R$ on $1-L$(up to sign). Thus desired result follows from similar arguments in Proposition~\ref{iso-moduli}. 
\end{enumerate}

\end{proof}

\subsection{Hilbert schemes of conics as Bridgeland moduli spaces}
In this subsection, we show that there is an isomorphism
\[p:C(X)\stackrel{\cong}{\to} \mathcal{M}_{\sigma'}(\mathcal{A}_{X}, 1-2L).\]
for $X:=X_{14}, X_{18}$ or $X_{22}$ and every Serre-invariant stability condition $\sigma'$ on $\mathcal{A}_X$. 

Let $\sigma'$ be a Serre-invariant Bridgeland stability condition on $\mathcal{A}_X$. 
First, we construct a natural morphism $p:\mathcal{C}(X)\to \mathcal{M}_{\sigma'}(\mathcal{A}_{X}, 1-2L)$:

\begin{lemma} \label{lem-4-1}
	The projection functor $\mathrm{pr}$ induces a morphism 
	$$p:\mathcal{C}(X)\to \mathcal{M}_{\sigma'}(\mathcal{A}_{X}, 1-2L)$$ where $p([C])=[\mathrm{pr}(I_C)]=[I_C]$ on the level of closed points.
\end{lemma} 

\begin{proof}
	From \cite[Lemma B.3.3]{KPS} we know that  $I_C\in \mathcal{A}_X$ for every conic $C$ on $X$. Hence $\mathrm{pr}(I_C)=I_C$. On the other hand, we know that $\mathrm{ext}^1(I_C, I_C)=2$, hence by Lemma \ref{ku-stable} $I_C$ is in the heart $\mathcal{A}_{\sigma'}$ up to some shifts and $\sigma'$-stable. Similar arguments in \cite[Section 5.1]{BMMS} and \cite[Section 5.1]{APR} shows that this shift can be chosen uniformly.

	Now we are going to show that there's a natural morphism $p$ induced by functor $\mathrm{pr}$. Since the projection functor $\mathrm{pr}: \mathrm{D^b}(X)\to \mathcal{A}_{X}$ is of Fourier-Mukai type, we can assume $K\in \mathrm{D^b}(X\times X)$ is the integral kernel and $\psi_K:\mathrm{D^b}(X)\xrightarrow{\mathrm{pr}} \mathcal{A}_{X}\hookrightarrow \mathrm{D^b}(X)$ be the Fourier-Mukai transform defined by $K$. Let  $\mathcal{I}$ be the universal ideal sheaf on $\mathcal{C}(X)\times X$. We define
	\[\psi':=\psi_K\times id_{\mathcal{C}(X)}=\psi_{K\boxtimes \mathcal{O}_{\Delta_{\mathcal{C}(X)}}}: \mathrm{D^b}(X\times \mathcal{C}(X))\to \mathrm{D^b}(X\times \mathcal{C}(X))\]
	Then $\psi'(\mathcal{I})$ is a family of objects in $\mathcal{A}_X$ parametrized by $\mathcal{C}(X)$, which defines a morphism $p:\mathcal{C}(X)\to \mathcal{M}_{\sigma'}(\mathcal{A}_{X}, 1-2L)$.
	
	To show $p([C])$ is given by $\mathrm{pr}(I_C)$ for any closed point $c=[C]\in \mathcal{C}(X)$, we denote $i_c: \{c\}\times X\to \mathcal{C}(X)\times X$. Then we have:
	$\psi_K(i_c^*(\mathcal{I}))\cong \mathrm{pr}(I_C)$
	and 
	\[i^*_c(\psi_{K\boxtimes \mathcal{O}_{\Delta_{\mathcal{C}(X)}}}(\mathcal{I}))\cong \psi_{i^*_cK\boxtimes \mathcal{O}_{\Delta_{\mathcal{C}(X)}}}(I_C)\cong \psi_K(I_C)=\mathrm{pr}(I_C)\]
	This means $p([C])$ is given by $[\mathrm{pr}(I_C)]$.
\end{proof}

Next, we show that $p$ is an isomorphism.

\begin{proposition} \label{conic-iso}
The morphism $p:C(X) \to  \mathcal{M}_{\sigma'}(\mathcal{A}_{X}, 1-2L)$ defined in Lemma \ref{lem-4-1} is an isomorphism.
\end{proposition}

\begin{proof}
	
	It's clear that $p$ is injective. The tangent map $\mathrm{d}p$ at a closed point $[C]$ is given by 
	\[\mathrm{d}p_{[C]}: \mathrm{Ext}^1(I_C, I_C)\to \mathrm{Ext}^1(\mathrm{pr}(I_C), \mathrm{pr}(I_C))\]
	Since $I_C$ is already in $\mathcal{A}_X$, the projection functor effects nothing on $I_C$. Hence we know that $\mathrm{d}p$ is an isomorphism at every closed point, and therefore $p$ is etale.
	
	Now since $\mathcal{C}(X)$ is projective and $\mathcal{M}_{\sigma'}(\mathcal{A}_{X}, 1-2L)$ is proper, $p$ is projective. Thus $p$ is an embedding. But $\mathcal{M}_{\sigma'}(\mathcal{A}_{X}, 1-2L)$ is irreducible and smooth by Proposition \ref{iso-moduli}, hence $p$ is actually an isomorphism.
	
\end{proof}

\begin{corollary} \label{406}
If $F\in \mathcal{A}_X$ is a $\sigma'$-stable object for some Serre-invariant stability conditions $\sigma'$ on $\mathcal{A}_X$ with $[F]=[I_C]\in \mathcal{N}(\mathcal{A}_X)$, then $F\cong I_C[2k]$ for some conics $C$ on $X$ and $k\in \mathbb{Z}$.
\end{corollary}


\section{Bridgeland moduli spaces of $(-1)$-classes on $X_{14}$ and $Y_3$} \label{sec-5}

In this section we fix $X:=X_{14}$ and $Y:=Y_3$. We show that $M_X(2,1,6)\cong \mathcal{M}_{\sigma'}(\mathcal{A}_{X_{14}}, 3s-t)$ and $M_X(2,-1,6)\cong \mathcal{M}_{\sigma'}(\mathcal{A}_{X_{14}}, 2s-t)$ on $X$.  And we have $M^b_Y(2,1,2)\cong \mathcal{M}_{\sigma}(\mathcal{K}u(Y), v-w)$ and $M^b_Y(2,-1,2)\cong \mathcal{M}_{\sigma}(\mathcal{K}u(Y), 2v-w)$ on $Y$.


\subsection{$M_X(2,1,6)$ as Bridgeland moduli space}

As shown in \cite{IM05}, there're two classes of sheaves in $M_X(2,1,6)$: globally generated locally free sheaves and non-locally free sheaves.

We first deal with locally free sheaves. When $E\in M_X(2,1,6)$ is a locally free sheaf, by \cite[Section 5.2]{IM05} we have an exact sequence:
$0\to \mathcal{O}_X\to E\to I_C(1)\to 0$
where $C$ is an elliptic sexic. 

First, we determine the image of $E$ under the projection functor $\mathrm{pr}$.

\begin{lemma} \label{vanish-bundle-1} Let $E\in M_X(2,1,6)$ be a locally free sheaf, then we have:
	
	\begin{enumerate}
		
		\item $h^0(E)=5, h^1(E)=h^2(E)=h^3(E)=0$.
		
		\item $h^0(E^{\vee})=h^1(E^{\vee})=h^2(E^{\vee})=h^3(E^{\vee})=0$.
		
		\item $\mathrm{hom}(E,E)=1, \mathrm{ext}^1(E,E)=2, \mathrm{ext}^2(E,E)=\mathrm{ext}^3(E,E)=0$.
	\end{enumerate}
\end{lemma}

\begin{proof}
	
	$(1)$ is from \cite[Lemma 5.1, Proposition 5.2]{IM05} and the fact $\chi(E)=5$.
	
	For $(2)$, by Serre duality we have $h^1(E^{\vee})=h^2(E(-1))=h^2(E^{\vee})$, hence $h^1(E^{\vee})=h^2(E^{\vee})=0$ by \cite[Proposition 5.2]{IM05}. And $h^0(E^{\vee})=0$ from stability of $E$, therefore by Serre duality we have $h^0(E^{\vee})=h^3(E(-1))=h^3(E^{\vee})=0$.
	
	$(3)$ is from \cite[Proposition 5.10]{IM05}.
\end{proof}

\begin{lemma} \label{eee}
Let $E\in M_X(2,1,6)$ be a locally free sheaf, then we have $\mathrm{ext}^i(\mathcal{E}^{\vee}, E)=0$ for all $i$.
\end{lemma}

\begin{proof}
 From the stability and Serre duality, we know $\mathrm{ext}^i(\mathcal{E}^{\vee}, E)=0$ for $i=0, 3$.    Since $E$ is globally generated, a general section of $E$ will vanish along an elliptic sexitc $C$ in $\mathrm{Gr}(2,6)$. Then we have an exact sequence as in \cite{IM05}:
	\[0\to \mathcal{O}_X\to E\to I_C(1)\to 0\]
	Applying Hom$(\mathcal{E}^{\vee},-)$ to this sequence, 
	and since $\mathrm{ext}^i(\mathcal{E}^{\vee}, \mathcal{O}_X)=0$ for every $i$, we have \[\text{Ext}^i(\mathcal{E}^{\vee}, E)=\text{Ext}^i(\mathcal{E}^{\vee}, I_C(1))=H^i(X, \mathcal{E}^{\vee}\otimes I_C), \forall i\geq 0.\]
	As in \cite[Section 5.2]{IM05}, $\mathcal{E}^{\vee}|_C$ is split of type $(3,3)$ or unsplit. When $\mathcal{E}^{\vee}|_C$ is split, it's clear that $h^0(\mathcal{E}^{\vee}|_C)=6$. When $\mathcal{E}^{\vee}|_C$ is unsplit, in this case $\mathcal{E}^{\vee}|_C$ is the tensor of a degree 3 line bundle $N$ with a unique vector bundle $F$ obtained as a non-trivial extension of $\mathcal{O}_C$. Thus $\mathcal{E}^{\vee}|_C$ is a non-trivial extension of $N$, which gives $h^0(\mathcal{E}^{\vee}|_C)=6$. Now if the restriction map $H^0(\mathcal{E}^{\vee})\to H^0(\mathcal{E}^{\vee}|_C)$ is not injective, then $C$ is contained in a copy of $\mathrm{Gr}(2,5)$. But this is impossible since $C$ is of degree $6$ and $\mathrm{Gr}(2,5)$ has degree $5$. Then the restriction map $H^0(\mathcal{E}^{\vee})\to H^0(\mathcal{E}^{\vee}|_C)$ is an isomorphism from the injectivity. Now applying $\Gamma(X, \mathcal{E}^{\vee}\otimes-)$ to the standard exact sequence of $C$, we obtain $\text{Ext}^i(\mathcal{E}^{\vee}, E)=H^i(\mathcal{E}^{\vee}\otimes I_C)=0$.
\end{proof}

Therefore by Lemma \ref{eee}, we have $\mathrm{pr}(E)=L_{\mathcal{O}_X}(E)$. From Lemma \ref{vanish-bundle-1} we know that $\mathrm{pr}(E)=L_{\mathcal{O}_X}(E)$ is given by $\ker(ev)[1]$, where $ev$ is the evaluation map:
\begin{equation} \label{ev-1}
	0\to \ker(ev)\to H^0(E)\otimes \mathcal{O}_X\xrightarrow{ev} E\to 0
\end{equation}
It's clear that $\mathrm{ch}(\ker(ev))=3-H-L+\frac{2}{3}P$. We are going to check that $\mathrm{pr}(E)=\ker(ev)[1]$ is $\sigma'$-stable.

\begin{lemma} \label{216-stable}
	We have $\mathrm{ext}^1(\ker(ev),\ker(ev))=2$. Hence $\mathrm{pr}(E)=\ker(ev)[1]\in \mathcal{A}_X$ is $\sigma'$-stable with respect to every Serre-invariant stability condition $\sigma'$ on $\mathcal{A}_X$.
\end{lemma}

\begin{proof}

	The second statement follows from the first one and Lemma \ref{ku-stable}, hence we only need to show the first statement.

Tensoring $\mathcal{O}_X(-1)$ with the sequence \ref{ev-1} and taking cohomology, we have $h^i(\ker(ev)(-1))=5\cdot h^i(\mathcal{O}(-1))$ for all $i$ by Lemma \ref{vanish-bundle-1}. Hence we know  $h^3(\ker(ev)(-1))=5$ and $h^i(\ker(ev)(-1))=0$ for $i\neq 3$. Thus by Serre duality we have $\mathrm{hom}(\ker(ev), H^0(E)\otimes \mathcal{O}_X)=25, \mathrm{ext}^i(\ker(ev), H^0(E)\otimes \mathcal{O}_X)=0$ for $i\neq 0$.

By Serre duality we have $\mathrm{Ext}^i(\ker(ev), E)=\mathrm{Ext}^{3-i}(E(1), \ker(ev))$.
Again by Serre duality, we have $\mathrm{ext}^i(E(1), H^0(E)\otimes \mathcal{O}_X)=0$ when $i\neq 3$ and  $\mathrm{ext}^3(E(1), H^0(E)\otimes \mathcal{O}_X)=25$. Also we have $\mathrm{Ext}^i(E(1), E)=\mathrm{Ext}^{3-i}(E, E)$, hence by Lemma \ref{vanish-bundle-1} we know that  $\mathrm{hom}(E(1), E)=\mathrm{ext}^1(E(1), E)=0$. Therefore applying $\mathrm{Hom}(E(1),-)$ to the exact sequence \ref{ev-1} and taking long exact sequence, we obtain $\mathrm{hom}(E(1),\ker(ev))=\mathrm{ext}^1(E(1), \ker(ev))=\mathrm{ext}^2(E(1), \ker(ev))=0$.
	
Finally, applying $\mathrm{Hom}(\ker(ev),-)$ to the sequence \ref{ev-1}, we obtain a long exact sequence.
From computations above, we have $\mathrm{ext}^i(\ker(ev), \ker(ev))=0$ for $i=2,3$. Since $\chi(\ker(ev), \ker(ev))=-1$, we only need to show $\mathrm{hom}(\ker(ev), \ker(ev))=1$. Note that $H^0(X, \ker(ev))=0$, then via the inclusion $\bigwedge^p \ker(ev)\hookrightarrow \bigwedge^{p-1} \ker(ev)\otimes H^0(E)$ we know that $H^0(X, \bigwedge^p \ker(ev))=0$ for all $p\geq 1$. And we have  $-1<\mu(\bigwedge^p \ker(ev))=p\mu(\ker(ev))<0$ for every $1\leq p\leq 2$. Therefore $\ker(ev)$ is a stable bundle by Hoppe's criterion \cite{hop} (see also \cite[Theorem 1.2]{Hoppe}). Hence we obtain that $\mathrm{hom}(\ker(ev), \ker(ev))=1$ and $\mathrm{ext}^1(\ker(ev), \ker(ev))=2$.
\end{proof}

\[\]

When $E\in M_X(2,1,6)$ is not locally free, by \cite[Proposition 5.11]{IM05}, there's an exact sequence:
\begin{equation} \label{511}
	0\to E\to T\to \mathcal{O}_L\to 0
\end{equation}
where $L$ is a line on $X$. 
By \cite[Proposition 3.5]{BCDvectorbundlesongenus7}, we know $T$ is a rank 2 globally generated stable vector bundle with $c_1(T)=H, c_2(T)=5L$. Thus from Proposition \ref{215}, we have $T\cong \mathcal{E}^{\vee}$.

\begin{lemma} \label{vanish-nonbundle} Let $E\in M_X(2,1,6)$ be a non-locally free sheaf, then we have:
	
	\begin{enumerate}
		
		\item $h^0(E)=5, h^1(E)=h^2(E)=h^3(E)=0$.
		
		\item $h^0(E(-1))=h^1(E(-1))=h^2(E(-1))=h^3(E(-1))=0$.
		
		\item $\mathrm{hom}(E,E)=1, \mathrm{ext}^1(E,E)=2, \mathrm{ext}^2(E,E)=\mathrm{ext}^3(E,E)=0$.
	\end{enumerate}
\end{lemma}

\begin{proof}
	(1) and (2) both follow from Proposition \ref{215} and the above exact sequence \ref{511}. (3) is from \cite[Proposition 5.12]{IM05}.
\end{proof}

\begin{lemma} \label{54}
Let $E\in M_X(2,1,6)$ be a non-locally free sheaf, then we have	$\mathrm{Ext}^i(\mathcal{E}^{\vee}, E)=0$ for every $i$.
\end{lemma}

\begin{proof}
	By Proposition \ref{215}, we know $\mathrm{ext}^1(T,T)=\mathrm{ext}^2(T,T)=\mathrm{ext}^3(T,T)=0$ and $\mathrm{hom}(T,T)=1$. From \cite[Proposition 5.11]{IM05} we have $T^{\vee}\otimes \mathcal{O}_L\cong \mathcal{O}_L\oplus \mathcal{O}_L(-1)$, hence we obtain $\mathrm{hom}(T, \mathcal{O}_L)=1$ and  $\mathrm{ext}^1(T, \mathcal{O}_L)=0$. Since $\mathrm{Hom}(T,T)\to \mathrm{Hom}(T, \mathcal{O}_L)$ is induced by the surjective map $T\twoheadrightarrow \mathcal{O}_L$, this map is non-trivial, hence for dimensional reason it's an isomorphism. Thus applying $\mathrm{Hom}(T,-)$ to the sequence \ref{511} above, from the long exact sequence we obtain that $\mathrm{Ext}^i(T,E)=0, \forall i\geq 0$.
	
\end{proof}

Therefore by Lemma \ref{54}, $\mathrm{pr}(E)=L_{\mathcal{O}_X}(E)$, which is given by $\ker(ev)[1]$, where $ev$ is the evaluation map:
\[0\to \ker(ev)\to H^0(E)\otimes \mathcal{O}_X\xrightarrow{ev} E\to 0\]
Now we want to check $\mathrm{pr}(E)=\ker(ev)[1]$ is $\sigma'$-stable. This is almost the same as locally free case: 

\begin{lemma} \label{216 stable nonfree}
	We have $\mathrm{ext}^1(\ker(ev),\ker(ev))=2$. Hence $\mathrm{pr}(E)=\ker(ev)[1]\in \mathcal{A}_X$ is $\sigma'$-stable with respect to every Serre-invariant stability condition $\sigma'$ on $\mathcal{A}_X$. 
\end{lemma}

\begin{proof}
The second statement follows from Lemma \ref{ku-stable}, hence we only need to show the first statement. Using Lemma \ref{vanish-nonbundle}, the same arguments in Lemma \ref{216-stable} shows that $\mathrm{ext}^i(\ker(ev), \ker(ev))=0$ for $i=2,3$. Since $\chi(\ker(ev), \ker(ev))=-1$, we only need to show $\mathrm{hom}(\ker(ev), \ker(ev))=1$.

Note that in this case $\ker(ev)$ is reflexive, hence is determined by the complement of any closed subset of codimension $\geq 2$ (see for example, \cite{har80}). And since the non-locally free locus of $\ker(ev)$ has codimension $\geq 3$, without loss of generality we can assume $\ker(ev)$ is locally free. Hence the same argument in Lemma \ref{216-stable} shows that $\mathrm{hom}(\ker(ev), \ker(ev))=1$.
\end{proof}

\begin{lemma} \label{216-induce-mor}
	The projection functor induces a morphism $$q:M_X(2,1,6)\to \mathcal{M}_{\sigma'} (\mathcal{A}_{X}, 3s-t)$$
	 for every Serre-invariant stability condition $\sigma'$ on $\mathcal{A}_X$, such that $q([E])=[\mathrm{pr}(E)]$.
\end{lemma}

\begin{proof}
	Since $\chi(E,E)=-1$, by \cite[Theorem 4.6.5]{HL} we know $M_X(2,1,6)$ is a fine moduli space. Thus using Lemma \ref{216-stable} and Lemma \ref{216 stable nonfree}, this lemma follows from the same argument in Lemma \ref{lem-4-1}.
	
	
\end{proof}


\begin{proposition}  \label{216-iso}
The morphism $q:M_X(2,1,6)\to \mathcal{M}_{\sigma'} (\mathcal{A}_{X}, 3s-t)$  defined in Lemma \ref{216-induce-mor} is an isomorphism.
\end{proposition}

\begin{proof}
	
	It's clear that $q$ is injective. Let $[E]\in M_X(2,1,6)$ be a closed point, then $T_{[E]}M_X(2,1,6)=\mathrm{Ext}^1(E,E)$. And we know $T_{[q(E)]}\mathcal{M}_{\sigma'}=\mathrm{Ext}^1(\ker(ev),\ker(ev))=\mathrm{Ext}^1(\ker(ev)[1], \ker(ev)[1])=\mathrm{Ext}^1(\mathrm{pr}(E), \mathrm{pr}(E))$. By definition 
	of projection, we have an exact triangle
	$$H^0(E)\otimes \mathcal{O}_X\xrightarrow{ev} E\to \mathrm{pr}(E)$$
	Applying $\mathrm{Hom}(E, -)$ to this triangle and using Lemma \ref{vanish-nonbundle}, we obtain a long exact sequence:
	\[0\to \mathrm{Hom}(E,E)\to \mathrm{Hom}(E, \mathrm{pr}(E))\to 0\to \mathrm{Ext}^1(E,E)\to \mathrm{Ext}^1(E, \mathrm{pr}(E))\to 0\]
	Since $(\mathrm{pr}, i)$ is an adjoint pair, where $i$ is an embedding of admissible subcategory. Thus the tangent map $\mathrm{d}_{[E]}q: \mathrm{Ext}^1(E, E)\to \mathrm{Ext}^1(E, \mathrm{pr}(E))=\mathrm{Ext}^1(\mathrm{pr}(E), \mathrm{pr}(E))$
	is an isomorphism.
	
	Finally, we know that $M_X(2,1,6)$ and $\mathcal{M}_{\sigma'} (\mathcal{A}_{X}, 3s-t)$ are both proper, hence $q$ is also proper. And since $q$ is injective and etale, it is an embedding. By Proposition \ref{iso-moduli} and Corollary \ref{all moduli spaces of (-1)-class isom}, we know that $\mathcal{M}_{\sigma'} (\mathcal{A}_{X}, 3s-t)$ is irreducible and smooth, therefore $q$ is actually an isomorphism.
	
\end{proof}

\subsection{$M_X(2,-1,6)$ as Bridgeland moduli space}



\begin{lemma} \label{2-16 in ku} Let $E\in M_X(2,-1,6)$, then we have:
	
	\begin{enumerate}
		\item $H^i(E)=0, \forall i\geq 0$.
		
		\item $\mathrm{Ext}^i(\mathcal{E}^{\vee}, E)=0, \forall i\geq 0$.
	\end{enumerate}
\end{lemma}

\begin{proof}
(1): Since $E(1)\in M_X(2,1,6)$, this follows from Lemma \ref{vanish-bundle-1} and Lemma \ref{vanish-nonbundle}.
	
	(2): 
	When $E$ is locally free, we have  $\mathrm{Ext}^i(\mathcal{E}^{\vee}, E)=\mathrm{Ext}^{3-i}(\mathcal{E}^{\vee}, E^{\vee})$. Then the statement follows from Lemma \ref{eee}. 
	
	When $E$ is not locally free, by Maruyama's restriction theorem we can choose a sufficiently general linear section $S$, such that $E(1)|_S$ and $\mathcal{E}^{\vee}|_S$ are both $\mu$-semistable. Since these two sheaves correspond to primitive vectors, $E(1)|_S$ and $\mathcal{E}^{\vee}|_S$ are actually $\mu$-stable. 
	From Lemma \ref{54} we know $\mathrm{ext}^i(\mathcal{E}^{\vee}, E(1))=0$ for all $i$. Since $\mathcal{E}^{\vee}|_S$ and $E(1)|_S$ are stable and not isomorphic to each other, we have $\mathrm{Hom}(\mathcal{E}^{\vee}, E(1)|_S)=0$ and $\mathrm{Ext}^2(\mathcal{E}^{\vee}, E(1)|_S)=\mathrm{Ext}^2(\mathcal{E}^{\vee}|_S, E(1)|_S)=\mathrm{Hom}(E(1)|_S, \mathcal{E}^{\vee}|_S)=0$.
	Now since $\chi(\mathcal{E}^{\vee}, E(1)|_S)=0$, we have  $\mathrm{Ext}^1(\mathcal{E}^{\vee}, E(1)|_S)=0$. Thus applying $\mathrm{\mathrm{Hom}}(\mathcal{E}^{\vee},-)$ to the short exact sequence $0\to E\to E(1)\to E(1)|_S\to 0$, we obtain $\mathrm{Ext}^i(\mathcal{E}^{\vee}, E)=0, \forall i\geq 0$.
\end{proof}

Therefore we know that $E\in \mathcal{A}_{X}$, which means $\mathrm{pr}(E)=E$.

\begin{lemma} \label{2-16 stable}
	Let $E\in M_X(2,-1,6)$, then $\mathrm{hom}(E,E)=1, \mathrm{ext}^1(E,E)=2$ and  $\mathrm{ext}^i(E,E)=0$ for $i\geq 2$. Hence $E$ is $\sigma'$-stable for every Serre-invariant stability condition $\sigma'$ on $\mathcal{A}_X$.
\end{lemma}

\begin{proof}
The second statement follows from Lemma \ref{ku-stable} and the first one. For the first statement, we have $\mathrm{Ext}^i(E,E)=\mathrm{Ext}^i(E(-1),E(-1))$, hence the statement follows from Lemma \ref{vanish-bundle-1} and Lemma \ref{vanish-nonbundle}.
\end{proof}

\begin{proposition}  \label{2-16-iso}
The projection functor $\mathrm{pr}$ induces an isomorphism 
$$p':M_X(2,-1,6)\xrightarrow{\cong} \mathcal{M}_{\sigma'} (\mathcal{A}_{X}, 2s-t)$$
\end{proposition}

\begin{proof}
	Since $\chi(E, E)=-1$ for $E\in M_X(2,-1,6)$, using \cite[Theorem 4.6.5]{HL} we know that $M_X(2,-1,6)$ is a fine moduli space. Thus using Lemma \ref{2-16 in ku} and Lemma \ref{2-16 stable}, the same argument in Lemma \ref{lem-4-1} shows that the projection functor $\mathrm{pr}$ induces a morphism $p':M_X(2,-1,6)\to \mathcal{M}_{\sigma'} (\mathcal{A}_{X}, 2s-t)$ such that $p'([E])=[\mathrm{pr}(E)]$ on the level of closed points.
	
	Now since $p'([E])=[\mathrm{pr}(E)]=[E]$, $p'$ is injective.
	Since the projection functor $\mathrm{pr}$ effects nothing on $E$, $p'$ is etale.
	And we know $p'$ is proper from that moduli spaces on both sides are proper. Therefore, $p'$ is proper injective and etale, which is an embedding. By Theorem \ref{all moduli spaces of (-1)-class isom}, $\mathcal{M}_{\sigma'}(\mathcal{A}_X, 2s-t)$ is irreducible and smooth, thus $p'$ is an isomorphism.
\end{proof}

\subsection{$M^{b}_Y(2,1,2)$ as Bridgeland moduli space}

Recall that a vector bundle is called normalised if $h^0(E(-1))=0$ and $h^0(E)\neq 0$. We denote by $M^{b}_Y(2,1,2)$ the moduli space of (semi)stable bundle of rank 2 and Chern class $c_1=1, c_2=2, c_3=0$ on $Y$.

\begin{lemma}
Every $E\in M^b_Y(2,1,2)$ is normalised.
\end{lemma}

\begin{proof}
Since $E$ is stable, we have $h^0(E(-1))=\mathrm{hom}(\mathcal{O}(1), E)=0$. A similar argument in \cite[Lemma 5.1]{IM05} shows that $h^2(E)=h^3(E)=0$, hence $h^0(E)\neq 0$ follows from $\chi(E)=3$.
\end{proof}

Thus by \cite[Lemma 1]{212cubic}, every $E\in M^b_Y(2,1,2)$ is ACM and indecomposable.

\begin{lemma} \label{pr-cone}
For every $E\in M^b_Y(2,1,2)$, we have $h^*(E(-1))=0$ and $h^0(E)=3$, $h^i(E)=0$ for $i\neq 0$.
\end{lemma}

\begin{proof}
       Since $E$ is ACM, we have $h^i(E(-1))=h^i(E)=0$ for $i\neq 0,3$. And from stability we have $h^0(E(-1))=0$. By Serre duality and stability of $E$, we have $h^3(E)=h^0(E^{\vee}(-2))=h^0(E(-3))=0$. Since $\chi(E(-1))=0$ and $\chi(E)=3$, we know that  $h^3(E(-1))=0$ and $h^0(E)=3$.
\end{proof}

Therefore by Lemma \ref{pr-cone}, $\mathrm{pr}(E)$ is given by $\mathrm{cone}(ev)$, where $ev: H^0(E)\otimes \mathcal{O}_Y\to E$ is the evaluation map.

\begin{lemma}  \label{cone-stable}
We have $\mathrm{ext}^1(\mathrm{cone}(ev),\mathrm{cone}(ev))=2$. Hence $\mathrm{pr}(E)=\mathrm{cone}(ev)$ is $\sigma$-stable for every Serre-invariant stability condition $\sigma$ on $\mathcal{K}u(Y)$.
\end{lemma}

\begin{proof}
Using \cite[Lemma 3]{212cubic}, a similar argument in Lemma \ref{216-stable} shows that $\mathrm{ext}^2(\mathrm{cone}(ev),\mathrm{cone}(ev))=\mathrm{ext}^3(\mathrm{cone}(ev),\mathrm{cone}(ev))=0$. Thus the statement follows from $\mathrm{hom}(\mathrm{cone}(ev), \mathrm{cone}(ev))=1$. 

\end{proof}

\begin{proposition} \label{212-iso}
    The projection functor $\mathrm{pr}$ induces an isomorphism 
    \[M^b_Y(2,1,2)\xrightarrow{\cong} \mathcal{M}_{\sigma}(\mathcal{K}u(Y), v-w)\]
\end{proposition}

\begin{proof}
  Using Lemma \ref{pr-cone} and Lemma \ref{cone-stable}, a similar argument in Lemma \ref{216-induce-mor} shows that $\mathrm{pr}$ induces a morphism $M^b_Y(2,1,2)\to \mathcal{M}_{\sigma}(\mathcal{K}u(Y), v-w)$ such that maps $[E]$ to $[\mathrm{pr}(E)]$. And as in Proposition \ref{216-iso}, we know that this morphism is etale, injective and proper, hence is an embedding. Since $\mathcal{M}_{\sigma}(\mathcal{K}u(Y), v-w)$ is irreducible and smooth by Corollary \ref{all moduli spaces of (-1)-class isom}, this is an isomorphism.
\end{proof}

\subsection{$M^{b}_Y(2,-1,2)$ as Bridgeland moduli space}

We denote by $M^{b}_Y(2,-1,2)$ the moduli space of (semi)stable bundle of rank 2 and Chern class $c_1=-1, c_2=2, c_3=0$ on $Y$.

\begin{lemma} \label{2--1-2-in-ku}
For every $E\in M^b_Y(2,-1,2)$, we have $E\in \mathcal{K}u(Y)$.
\end{lemma}

\begin{proof}
Since $E(1)\in M^b_Y(2,1,2)$, by Lemma \ref{pr-cone} we have $h^*(E)=0$. And by Serre duality we have $h^i(E(-1))=h^{3-i}(E^{\vee}(-1))=h^{3-i}(E)$, thus from Lemma \ref{pr-cone} we know $h^*(E(-1))=0$.
\end{proof}

\begin{lemma} \label{2--1-2-stable}
For every $E\in M^b_Y(2,-1,2)$, we have $\mathrm{hom}(E, E)=1, \mathrm{ext}^1(E,E)=2, \mathrm{ext}^2(E,E)=\mathrm{ext}^3(E,E)=0$. Hence $\mathrm{pr}(E)=E$ is $\sigma$-stable with respect to  every Serre-invariant stability condition $\sigma$ on $\mathcal{K}u(Y)$.
\end{lemma}

\begin{proof}
Since $E(1)\in M^b_Y(2,1,2)$, this is from $\mathrm{hom}(E,E)=1$ and \cite[Lemma 3]{212cubic}.
\end{proof}

\begin{proposition} \label{2-12-iso}
    The projection functor $\mathrm{pr}$ induces an isomorphism 
    \[M^b_Y(2,-1,2)\xrightarrow{\cong} \mathcal{M}_{\sigma}(\mathcal{K}u(Y), 2v-w)\]
\end{proposition}

\begin{proof}
  Using Lemma \ref{2--1-2-in-ku} and Lemma \ref{2--1-2-stable}, a similar argument in Proposition \ref{2-16-iso} shows that $\mathrm{pr}$ induces a morphism $M^b_Y(2,-1,2)\to \mathcal{M}_{\sigma}(\mathcal{K}u(Y), 2v-w)$ such that maps $[E]$ to $[\mathrm{pr}(E)]=[E]$. And as in Proposition \ref{2-16-iso}, we know that this morphism is etale, injective and proper, hence is an embedding. Since $\mathcal{M}_{\sigma}(\mathcal{K}u(Y), 2v-w)$ is irreducible and smooth by Corollary \ref{all moduli spaces of (-1)-class isom}, this is an isomorphism.
\end{proof}

\section{Instanton sheaves and wall-crossing} \label{sec-6}

Let $Y:=Y_d$ be a prime Fano threefold of index 2 and degree $d$. The notion of mathematical instanton bundle was first introduced on $\mathbb{P}^3$ and generalized to $Y$ in \cite{fae11} and \cite{kuz12}.

\begin{definition}
    Let $Y$ be a prime Fano threefold of index 2. An \emph{instanton bundle of charge n} on $Y$ is a Gieseker-stable vector bundle $E$ of rank 2 with $c_1(E)=0$, $c_2(E)=n$, and satisfies the instantonic condition:
    \[H^1(Y, E(-1))=0\]
\end{definition}

As shown in \cite{D} and \cite{qin2019moduli}, every Gieseker-semistable sheaf with $(r,c_1,c_2,c_3)=(2,0,2,0)$ on $Y_3$ and $Y_4$ automatically satisfies $H^1(E(-1))=0$. Therefore, we can give a more general definition:

\begin{definition}
    Let $Y$ be a prime Fano threefold of index 2. An \emph{instanton sheaf} on $Y$ is a Gieseker-semistable sheaf $E$ of rank 2 with $c_1(E)=0$, $c_2(E)=2$ $c_3(E)=0$, and satisfies the instantonic condition:
    \[H^1(Y, E(-1))=0\]
\end{definition}

\subsection{Wall-crossing}

In this subsection, we show that when $d\neq 2$, the instantonic condition is equivalent to the non-existence of the maximal semicircle wall with respect to $\sigma_{\alpha, \beta}$.

\begin{lemma} \label{leq-2}
  Let $E\in \mathrm{Coh}^{-\frac{1}{2}}(Y)$ be a $\sigma_{\alpha, -\frac{1}{2}}$-semistable object for some $\alpha>0$ with $\mathrm{ch}_{\leq 2}(E)=(1,0,-\frac{1}{d}H^2)$. Then $E$ is $\sigma_{\alpha, -\frac{1}{2}}$-semistable for every $\alpha>0$ and $\mathrm{ch}(E)=1-L+xP$ for $x\in \mathbb{Z}_{\leq 0}$ when $d\neq 1$, and $x\in \mathbb{Z}_{\leq 1}$ when $d=1$.
\end{lemma}

\begin{proof}
       Assume $\mathrm{ch}(E)=1-L+xP$ for $x\in \mathbb{Z}$. 
       From the wall-crossing argument in \cite[Proposition 4.1]{PY20} we know that $E$ is $\sigma_{\alpha, -\frac{1}{2}}$-semistable for every $\alpha>0$. Hence by \cite{Li15}, \cite[Conjecture 4.1]{BMS} holds for $E$ and  $\alpha=0$. Therefore when $d\neq 1$, we have $x<1$, which means $x\leq 0$. When $d=1$, we have $x\leq \frac{3}{2}$, which means $x\leq 1$
\end{proof}

Next we rule out all possible walls of $[E]=2-2L$ on $\beta=-\frac{1}{2}$.

\begin{proposition} \label{wall-cross-1}
When $d\neq 5$, the only possible walls for $\mathrm{ch}(E)=2-2L$ on $\beta=-\frac{1}{2}$ are realized by extensions of ideal sheaves of lines. 

When $d=5$, there're two possible walls on $\beta=-\frac{1}{2}$, realized by extensions of ideal sheaves of lines and object $\mathcal{O}(-1)[1]$. 

\end{proposition}


	
\begin{proof}
	

	If there's a wall, then it is given by a short exact sequence
	\[0\to A\to E\to B\to 0\]
	such that:
	\begin{enumerate}
		\item $\mu_{\alpha, -\frac{1}{2}}(A)=\mu_{\alpha, -\frac{1}{2}}(B)=\mu_{\alpha, -\frac{1}{2}}(E)$;
		
		\item $\Delta_H(A)\geq 0, \Delta_H(B)\geq 0$;
		
		\item $\Delta_H(A)\leq \Delta_H(E), \Delta_H(B)\leq \Delta_H(E)$
		
		\item $\mathrm{ch}^{-\frac{1}{2}}(A)+\mathrm{ch}^{-\frac{1}{2}}(B)=\mathrm{ch}^{-\frac{1}{2}}(E)$;
	\end{enumerate}
	
	For (4) we can assume
	\[(2, H, \frac{d-8}{4d}H^2)=(a, \frac{b}{2}H, \frac{c}{8d}H^2)+ (2-a, \frac{2-b}{2}H, \frac{2d-16-c}{8d}H^2)\]
	for some $a,b,c\in \mathbb{Z}$ such that $2|b-a$   (this is from $\mathrm{ch}^{-\frac{1}{2}}_1-\frac{a}{2}H=\mathrm{ch}_1=c_1\in \mathbb{Z}\cdot H$). Since $A,B \in \mathrm{Coh}^{-\frac{1}{2}}(Y)$, we have Im$(Z(A))\geq 0$ and Im$(Z(B))\geq 0$. Thus $b\geq 0$ and $2-b\geq 0$, which implies $b=0$ or $b=1$ or $b=2$.
	
	Since $\mu_{\alpha, -\frac{1}{2}}(E)=\frac{d-8}{4d}-\alpha^2$, $\Delta_H(E)=\frac{8}{d}(H^3)^2$, the previous conditions can be written as:
	
	\begin{enumerate}
		\item $\frac{1}{b}(\frac{c}{4d}-\alpha^2 a)=\frac{d-8}{4d}-\alpha^2=\frac{1}{2-b}(\frac{2d-16-c}{4d}-\alpha^2(2-a))$;
		
		\item $\frac{b^2}{4}-\frac{ac}{4d}\geq 0, ~ \frac{(2-b)^2}{4}+\frac{(2-a)(c+16-2d)}{4d}\geq 0$;
		
		\item $\frac{b^2}{4}-\frac{ac}{4d}\leq \frac{8}{d}, ~ \frac{(2-b)^2}{4}+\frac{(2-a)(c+16-2d)}{4d}\leq \frac{8}{d}$
	\end{enumerate}
	
Now from $\frac{d-8}{4d}-\alpha^2<\infty$ for all $\alpha >0$, we know that $b\neq 0$ or $2$.
	Then $b=1$, and therefore $c-d+8=4d(a-1)\alpha^2$ by (1). Since $2|b-a$, we have $a$ is odd. By (2) and (3) we have $d-32\leq ac\leq d$ and $d-32\leq (a-2)(c+16-2d)\leq d$. 
	Now note that $\mathrm{ch}^{-\frac{1}{2}}_1(A)=\frac{1}{2}H=\mathrm{ch}_1(A)+\frac{1}{2}aH$, thus $\mathrm{ch}_1(A)=c_1(A)=\frac{1-a}{2}H$. And $\mathrm{ch}^{-\frac{1}{2}}_2(A)=\frac{c}{8d}H^2=\mathrm{ch}_2+\frac{1}{2}\mathrm{ch}_1H+\frac{1}{8}\mathrm{ch}_0H^2=\frac{c_1^2-2c_2}{2}+\frac{1-a}{4}H^2+\frac{1}{8}aH^2$, hence $(\frac{c}{8d}-\frac{(1-a)^2}{8}-\frac{2-a}{8})H^2=-c_2\in \mathbb{Z}\cdot H$. 
	In summary, we have following system of inequalities:
	
	\begin{enumerate}
	    \item $d=1,2,3,4,5$ and $a,c\in \mathbb{Z}$, $a$ odd and $\alpha>0$;
	    
	    \item $c-d+8=4d(a-1)\alpha^2$;
	    
	    \item $d-32\leq ac\leq d$;
	    
	    \item $d-32\leq (a-2)(c+16-2d)\leq d$;
	    
	    \item $\frac{c}{8}-\frac{d(1-a)^2}{8}-\frac{d(2-a)}{8}\in \mathbb{Z}$.
	\end{enumerate}
	
After solving this, we find that the only possible solutions are $(a,c, \alpha)=(1,d-8,\mathbb{R}_{>0})$ when $d\neq 5$, and $(a,c, \alpha)=(1.-3, \mathbb{R}_{>0}), (-1,-5, \frac{1}{\sqrt{20}}),(3,-1, \frac{1}{\sqrt{20}})$ when $d=5$.

	When $(a,c)=(1,d-8)$, the destabilizing sequence is given by 
	\[0\to A\to E\to B\to 0\]
	where $\mathrm{ch}_{\leq 2}^{-\frac{1}{2}}(A)=\mathrm{ch}_{\leq 2}^{-\frac{1}{2}}(B)=(1, \frac{1}{2}H, \frac{d-8}{8d}H^2)$. Since $\mathrm{ch}_3(A)+\mathrm{ch}_3(B)=\mathrm{ch}_3(E)=0$, by Lemma \ref{leq-2} we know actually $\mathrm{ch}(A)=\mathrm{ch}(B)=1-\frac{1}{d}H^2$ and $\sigma_{\alpha, -\frac{1}{2}}$-semistable for every $\alpha>0$. Hence from \cite[Lemma 2.7]{BMS} and $\mathrm{Pic}(Y)\cong \mathbb{Z}$, we know that $A$ and $B$ are $\mu$-semistable sheaves, hence both isomorphic to ideal sheaves of lines. 
	
	
	
	When $d=5$ and $(a,c)=(-1,-5),(3,-1)$, we have either $\mathrm{ch}^{-\frac{1}{2}}_{\leq 2}(A)=(-1,\frac{1}{2}H,-\frac{1}{8}H^2)$ or $\mathrm{ch}^{-\frac{1}{2}}_{\leq 2}(B)=(-1,\frac{1}{2}H,-\frac{1}{8}H^2)$. Thus by a standard result in \cite[Proposition 2.14]{BLMS} we know that $A\cong \mathcal{O}(-1)[1]$ or $B\cong \mathcal{O}(-1)[1]$.
\end{proof}

From a similar computation in Proposition \ref{wall-cross-1}, we have:

\begin{lemma} \label{no wall intersect -1}
Let $d\neq 2$ and $E\in \mathrm{Coh}^{-1}(Y)$ with $\mathrm{ch}(E)=2-2L$. Then there's no semicircle wall intersect with $\beta=-1$. Therefore, the maxiaml possible semicircle wall for $E$ is  $(\beta-\frac{d+2}{2d})^2+\alpha^2=(\frac{d-2}{2d})^2$ realized by $\mathcal{O}(-1)[1]$.
\end{lemma}




The following lemma gives another description of instantonic condition as existence of walls:

\begin{lemma} \label{inst-equiv-to-wall}
Let $d\neq 1, 2$ and $E\in M^{ss}_{Y}(2,0,2)$. Then the following are equivalent:

\begin{enumerate}
    \item $E$ satisfies instantonic condition.
    
    \item The wall given by $\mathcal{O}(-1)[1]$ with respect to $\sigma_{\alpha,\beta}$ is not an actual wall for $E$.
    
\end{enumerate}

When $d=1$, we have $(1)$ implies $(2)$.

\end{lemma}

\begin{proof}
Since $\chi(E(-1))=0$ and $h^0(E(-1))=h^3(E(-1))=0$ from stability, the instantonic condition is equivalent to $H^2(E(-1))=0$. By Serre duality, the instantonic condition is equivalent to $H^2(E(-1))=\mathrm{Hom}(\mathcal{O}(1), E[2])=\mathrm{Hom}(E, \mathcal{O}(-1)[1])=0$.


$(1)\Rightarrow (2)$: Since we have $\mathrm{Hom}(\mathcal{O}(-1)[1], E)=0$, these lead to the non-existence of the wall by definition.

$(2)\Rightarrow (1)$: The instantonic condition is equivalent to $H^2(E(-1))=\mathrm{Hom}(E, \mathcal{O}(-1)[1])=0$.  By \cite[Proposition 4.8]{bayer2020desingularization} we know $E$ is $\sigma_{\alpha, \beta}$-semistable for  $\alpha\gg0$ and $\beta<0$. Now by Lemma \ref{no wall intersect -1}, we know that $\mathcal{O}(-1)[1]$ gives the biggest semicircle wall. So if $\mathcal{O}(-1)[1]$ is not an actual wall, then after crossing this wall, $E$ is remain $\sigma_{\alpha, \beta}$-semistable. But when $d=3,4,5$, under this wall we have $E\in \mathrm{Coh}^{\beta}(Y)$ and $\mathcal{O}(-1)[1]\in \mathrm{Coh}^{\beta}(Y)$, and $\mu_{\alpha, \beta}(E)>\mu_{\alpha, \beta}(\mathcal{O}(-1)[1])$, which gives $\mathrm{Hom}(E,\mathcal{O}(-1)[1])=0$. 
\end{proof}


\begin{corollary} \label{inst-stable-first-tilt}
Let $E$ be an instanton sheaf. Then $E\in \mathrm{Coh}^{\beta}(Y)$ is $\sigma_{\alpha, \beta}$-semistable for every $\alpha>0$ and $\beta=-\frac{1}{2},-1$.
\end{corollary}

\begin{proof}
       This immediately follows from 
       Proposition \ref{wall-cross-1}, Lemma \ref{no wall intersect -1} and Lemma \ref{inst-equiv-to-wall}.
\end{proof}

We denote the moduli space of instanton sheaves by $M^{inst}_Y$. 
Let $Y=Y_3,Y_4$ or $Y_5$. We collect some properties and classifications of instanton sheaves from \cite{D}, \cite{kuz12}, \cite{qin2019moduli} and \cite{qinv5}.

\begin{lemma} \label{ins-prop} Let $Y=Y_3,Y_4$ or $Y_5$. Let $E\in M^{ss}_Y(2,0,2)$.
	\begin{enumerate} 
		\item When $d=3,4$, $E$ is an instanton sheaf.
		
		\item If $E$ is Gieseker-stable, then $E$ is either locally free or defined by
		\[0\to E\to H^0(\theta(1))\otimes \mathcal{O}_Y\to \theta(1)\to 0\]
		where $\theta$ is a theta-characteristic of a smooth conic $C\subset Y_3$. If $E$ is strictly Gieseker-semistable, then $E$ is an extension of two ideal sheaves of lines.
		
		
		\item For every $E\in M^{inst}_Y$, we have $E\in \mathcal{K}u(Y)$.
		
		\item $M^{inst}_Y$ is irreducible, projective and smooth of dimension $5$.
		
	\end{enumerate}
\end{lemma}

\begin{proof}
       (1): When $Y=Y_3$, this is proved in the proof of \cite[Theorem 2.4]{D}. When $Y=Y_4$, this is from \cite[Theorem 1.2]{qin2019moduli}.
       
       (2): This is from \cite[Theorem 3.5]{D}, \cite[Theorem 1.4]{qin2019moduli} and \cite[Theorem 1.2]{qinv5}.
       
       (3): This is from the classifications of such sheaves and the same argument as in \cite[Lemma 4.1]{qin2019moduli}.
       
       (4): This is from \cite[Theorem 4.6]{D}, \cite[Theorem 5.4]{qin2019moduli} and \cite[Theorem 5.6]{qinv5}.
\end{proof}

\section{Moduli of instanton sheaves on $Y_d$ as Bridgeland moduli space} \label{sec-7}

Let $Y:=Y_d$. In this section, we are going to show that the projection functor $\mathrm{pr}$ induces an isomorphism $M^{inst}_{Y}\xrightarrow{\cong} \mathcal{M}^{ss}_{\sigma}(\mathcal{K}u(Y), 2-2L)$ for $d\neq 1,2$ and an open immersion $M^{inst}_{Y}\hookrightarrow \mathcal{M}^{ss}_{\sigma}(\mathcal{K}u(Y), 2-2L)$ for $d=1,2$.

\subsection{In Kuznetsov Component}

First, we show that every instanton sheaf $E$ is in the Kuznetsov component $\mathcal{K}u(Y)$. When $d=3,4$ and $5$, this is shown  in \cite{D}, \cite{qin2019moduli} and \cite{qinv5} by using the classifications of instanton sheaves. We give another proof, which not need classification  results on $Y_1$ and $Y_2$, and also work for $d=3,4,5$.

\begin{lemma} \label{inst-in-ku}
	For every $E\in M^{inst}_Y$, we have $H^*(E)= H^*(E(-1))=0$. Thus $E\in \mathcal{K}u(Y)$.
\end{lemma}

	
	

\begin{proof}
       When $d= 3,4,5$, this is from Lemma \ref{ins-prop}.

       Assume $d=1,2$. From the stability and Serre duality, we have $h^0(E)=h^0(E(-1))=h^3(E)=h^3(E(-1))=0$. By the instantonic condition and $\chi(E(-1))$=0, we have $h^1(E(-1))=0$. Since $\chi(E)=0$, we only need to show $h^2(E)=0$.
       
       When $d=2$, from Corollary \ref{inst-stable-first-tilt} we know that $E\in \mathrm{Coh}^{-1}(Y)$ is $\sigma_{\alpha, -1}$-semistable for every $\alpha>0$. By \cite[Proposition 2.14]{BLMS} we have that $\mathcal{O}(-2)[1]\in \mathrm{Coh}^{-1}(Y)$ is $\sigma_{\alpha, -1}$-stable for every $\alpha>0$. 
       Thus taking $0< \alpha \ll 1$ gives $\mu_{\alpha, -1}(E)>\mu_{\alpha,-1}(\mathcal{O}(-2)[1])$, which implies $\mathrm{hom}(E, \mathcal{O}(-2)[1])=\mathrm{hom}(\mathcal{O}(-2)[1], E(-2)[3])=h^2(E)=0$.
       
       When $d=1$, from Lemma \ref{wall-cross-1}, Corollary \ref{inst-stable-first-tilt} and Lemma  \ref{inst-equiv-to-wall}, we know that $E\in \mathrm{Coh}^{-\frac{3}{2}}(Y)$ is $\sigma_{\alpha, -\frac{3}{2}}$-semistable for $\alpha\gg 0$, and remain semistable when $\alpha=\frac{1}{2}-\epsilon$ for $0<\epsilon\ll 1$. By \cite[Proposition 2.14]{BLMS} we have that $\mathcal{O}(-2)[1]\in \mathrm{Coh}^{-\frac{3}{2}}(Y)$ is $\sigma_{\alpha, -\frac{3}{2}}$-stable for every $\alpha>0$. 
       Thus taking $\alpha=\frac{1}{2}-\epsilon$ for $0<\epsilon\ll 1$ gives $\mu_{\alpha, -\frac{3}{2}}(E)>\mu_{\alpha,-\frac{3}{2}}(\mathcal{O}(-2)[1])$, which implies $\mathrm{hom}(E, \mathcal{O}(-2)[1])=\mathrm{hom}(\mathcal{O}(-2)[1], E(-2)[3])=h^2(E)=0$.
\end{proof}

\begin{proposition} \label{sheaf-stable}
	Let $E\in M^{inst}_Y$ be an instanton sheaf. Then we have $E[1]\in \mathcal{A}(\alpha, -\frac{1}{2})$ is $\sigma(\alpha, -\frac{1}{2})$-semistable for every $0<\alpha<\frac{1}{2}$. If $d\neq 1$ and $E$ is Gieseker-stable, then $E[1]$ is $\sigma(\alpha, -\frac{1}{2})$-stable for every $0<\alpha<\frac{1}{2}$.
\end{proposition}

\begin{proof}
       This is from Lemma \ref{inst-in-ku} and Corollary \ref{inst-stable-first-tilt}. When $d\neq 1$ and  $E[1]$ is strictly $\sigma(\alpha, -\frac{1}{2})$-semistable, by Lemma \ref{same-slope}, Lemma \ref{non-negative} and \cite[Theorem 1.1]{PY20} we know that $E$ is an extension of ideal sheaf of lines, hence is strictly Gieseker-semistable.
\end{proof}



\subsection{Bridgeland (semi)stable objects are (semi)stable sheaves}

In this subsection we assume $3\leq d\leq 5$.
We show that for every $\sigma(\alpha, -\frac{1}{2})$-(semi)stable object $F\in \mathcal{A}(\alpha, -\frac{1}{2})$ with $[F]=2L-2$, we have  $F[-1]$ is Gieseker-(semi)stable sheaf.






Our argument is as follows: 

First, using Proposition \ref{wall-cross-1}, a similar argument in \cite[Lemma 6.11]{BLMS} shows that when starting at $(\alpha_0, \beta_0)=(\frac{d-2}{2d}, -\frac{d+2}{2d})$ and $\beta_0$ approaching to $-\frac{1}{2}$, the only wall $\mathcal{C}$ will meet is realized by $\mathcal{O}(-1)^{\oplus a}[2]$ for some $a\in \mathbb{Z}_{>0}$. Then following the argument in \cite[Proposition 4.6]{PY20}, if $\mathcal{C}$ is not an actual wall, we can show that $F[-1]\in \mathrm{Coh}^{-\frac{1}{2}}(Y)$ is $\sigma_{\alpha, -\frac{1}{2}}$-semistable for every $\alpha>0$. Then again by \cite[Lemma 2.7]{BMS}, $F[-1]$ is a $\mu$-semistable sheaf. Therefore, by Lemma \ref{mu-ss-g-ss}, Lemma \ref{ins-prop} and the standard slope-comparison argument, $F[-1]$ is actually an instanton sheaf.

When $\mathcal{C}$ is an actual wall, we also follow the argument in \cite[Proposition 4.6]{PY20} and show that crossing this wall will lead to a contradiction.

\begin{lemma} \label{infty-wall}
Let $(\alpha_0, \beta_0)=(\frac{d-2}{2d},-\frac{d+2}{2d})$ and $F\in \mathrm{Coh}^0_{\alpha_0, \beta_0}(Y)$ be a $\sigma^0_{\alpha_0, \beta_0}$-semistable object with $\mathrm{ch}(F)=2L-2$. Then when $\beta_0$ approaching $\beta=-\frac{1}{2}$, the only wall of $F$ with respect to $\sigma^0_{\alpha, \beta}$ is realized by $\mathcal{O}(-1)^{\oplus a}[2]$ for some $a\in \mathbb{Z}_{>0}$.
\end{lemma}

\begin{proof}
     As shown in Lemma \ref{no wall intersect -1}, there's no semicircle wall intersect with $\beta=-1$. Thus we can follow the argument in \cite[Lemma 6.11]{BLMS}. 
     
     We consider the location of possible walls for $F$. They are given as the intersection of lines through $v^2(F)$ with the interior of the negative cone $\Delta(-)<0$. By Lemma \ref{no wall intersect -1}, no such wall can be in the interior of the triangle with vertices $v^2(F), v^2(\mathcal{O}(-1))$ and $(0,0,1)$. A wall is given by a sequence
     \[0\to A\to F\to B\to 0\]
     such that when $(\alpha,\beta)$ on the wall, then $Z^0_{\alpha, \beta}(A)$ and $Z^0_{\alpha, \beta}(B)$ lie on the open line segment connecting 0 and $Z^0_{\alpha, \beta}(F)$. By continuity, this still holds at the end point $(\alpha, \beta)=(0,-1)$. Assume $\mathrm{ch}^{-1}_{\leq 2}(A)=(a,bH,\frac{c}{2d}H^2)$ for $a,b,c\in \mathbb{Z}$. Since $\mu^0_{\alpha_0, \beta_0}(F)=+\infty$, we have $\mathrm{Re}(Z^0_{\alpha_0, \beta_0}(A))\leq 0$ and $\mathrm{Re}(Z^0_{\alpha_0, \beta_0}(A))\leq 0$. Thus we have $b=-2,-1,0$. First we assume $b=-1$, then we have  $\mu_{0,-1}(A)=-\frac{c}{2d}$. Since $\mathrm{ch}^{-1}_{\leq 2}(\mathcal{O}(-1))=(1,0,0)$, we have $\mu^0_{0,-1}(\mathcal{O}(-1))=+\infty$. Therefore, $\mu^0_{0,-1}(A)=\mu^0_{0,-1}(\mathcal{O}(-1))$ implies $c=0$.  Hence we have
     $\mathrm{ch}^{-1}_{\leq 2}(A)=(a,-H,0)$ and $\mathrm{ch}^{-1}_{\leq 2}(B)=(-2-a,-H,(\frac{2}{d}-1)H^2)$. But this is impossible since we must have  $\mu^0_{0,-1}(B)=+\infty$.
     Thus we have $b=0$ or $2$, and in these cases either $\mathrm{ch}_{\leq 2}(A)=(a,-aH,\frac{a}{2}H^2)$ or  $\mathrm{ch}_{\leq 2}(B)=(a,-aH,\frac{a}{2}H^2)$.
     By standard arguments, we conclude that $A\cong \mathcal{O}(-1)^{\oplus a}[2]$ or $B\cong \mathcal{O}(-1)^{\oplus a}[2]$.
\end{proof}

\begin{lemma} \label{mu-ss-g-ss}
Let $A$ be a $\mu$-semistable sheaf with $\mathrm{ch}(A)=2-2L$, such that $\sigma_{\alpha, -\frac{1}{2}}$-semistable for $\alpha \gg 0$, then $A$ is Gieseker-semistable.
\end{lemma}

\begin{proof}
If $A$ is not Gieseker-semistable, let $G\subset A$ be the destabilizing sheaf. Then $\mathrm{rk}(G)=1$, $\chi(G(n))>\chi(A(n))/2$ when $n\gg 0$ and $G$ is Gieseker-semistable. Assume that  $\mathrm{ch}(G)=1+aH+\frac{b}{2}L+\frac{c}{2}P$ where $a,b,c\in \mathbb{Z}$. 

	We know that $\chi(G(n))=(1+\frac{d+3}{3}a+\frac{b+c}{2})+(\frac{d+3}{3}+da+\frac{b}{2})n+(a+1)\frac{d}{2}n^2+\frac{d}{6}n^3$ and $\chi(A(n))/2=\frac{d}{3}n+\frac{d}{2}n^2+\frac{d}{6}n^3$. From $\chi(G(n))>\chi(A(n))/2$ for $n\gg 0$, we have:
	
	\begin{enumerate}
		\item $a>0$;
		
		\item $a=0$, $-2<b$;
		
		\item $a=0$, $b=-2$, $c>0$.
	\end{enumerate}
By the $\mu$-semistability of $A$, (1) is impossible.
	For (2) we have $\mu_{\alpha, -\frac{1}{2}}(G)>\mu_{\alpha, -\frac{1}{2}}(A)$, which contradicts with the $\sigma_{\alpha, -\frac{1}{2}}$-semistability of $A$. Hence the only possible case is $\mathrm{ch}(G)=1-L+\frac{c}{2}P$  and $c>0$. Since $\chi(G)\in \mathbb{Z}$, we can assume that $\mathrm{ch}(G)=1-L+xP$ for $x\in \mathbb{Z}_+$.  

Since $1-L+xP$ is a primitive class, $G$ is actually $\mu$-stable. Then by \cite[Lemma 2.7]{BMS}, $G$ is $\sigma_{\alpha, -\frac{1}{2}}$-stable for $\alpha \gg 0$. But by Lemma \ref{leq-2} this
contradicts with $x\in \mathbb{Z}_+$. Hence we conclude that $A$ is Gieseker-semistable.


	
	

\end{proof}

\begin{proposition} \label{stable-implies-sheaf}
Assume $d\neq 1,2$.	For every object $F\in \mathcal{A}(\alpha, -\frac{1}{2})$ with $[F]=2L-2$ and $\sigma(\alpha, -\frac{1}{2})$-(semi)stable for some $0<\alpha< \frac{1}{2}$, we have $F[-1]$ is a  (semi)stable instanton sheaf.
\end{proposition}

\begin{proof}
       We argue as \cite[Proposition 4.6]{PY20}. By \cite[Proposition 3.6]{PY20}, without loss of generality we can assume $F$ is in the heart $\mathrm{Coh}^0_{\alpha, \beta}(Y)$ for $(\alpha, \beta)\in V$.

       From $d\geq 3$ we can find a point $(\alpha_0, \beta_0)=(\frac{d-2}{2d},-\frac{d+2}{2d})\in V$ such that $F\in \mathrm{Coh}^0_{\alpha_0, \beta_0}(Y)$ and $\mu^0_{\alpha_0, \beta_0}=+\infty$. By Lemma \ref{no wall intersect -1} and Lemma \ref{infty-wall}, this is the top point of the semicircle wall, denote by $\mathcal{C}$, realized by $\mathcal{O}(-1)^{\oplus a}[2]$.
       
       Assume that $\mathcal{C}$ is not an actual wall for $F$. Thus $F$ is $\sigma^0_{\alpha_0, \beta_0}$-semistable and remains semistable when $\beta_0$ approaching $-\frac{1}{2}$.
       
     By the definition of $\mathrm{Coh}^0_{\alpha, -\frac{1}{2}}(Y)$, there's an exact triangle
	\[A[1]\to F\to B\]
	such that $A$ (resp. $B$) $\in \mathrm{Coh}^{-\frac{1}{2}}(Y)$ with $\sigma_{\alpha, -\frac{1}{2}}$-semistable factors having slope $\mu_{\alpha, -\frac{1}{2}}\leq 0$ (resp. $\mu_{\alpha, -\frac{1}{2}}>0$). Since $F$ is $\sigma^0_{\alpha, -\frac{1}{2}}$-semistable, we have that $Z_{\alpha, -\frac{1}{2}}(B)=0$. Hence  either $B$ support on a point or $B=0$, and therefore $\mathrm{ch}(A)=(2,0, -\frac{2}{d}H^2, mP)$ where $m\geq 0$ is the length of $B$.  Moreover, 
	$A[1]$ is $\sigma^0_{\alpha, -\frac{1}{2}}$-semistable and, since $A\in \mathrm{Coh}^{-\frac{1}{2}}(Y)$, 
	we have that $A$ is $\sigma_{\alpha, -\frac{1}{2}}$-semistable.
	
	But as we showed in Proposition \ref{wall-cross-1}, $A$ is actually $\sigma_{\alpha, -\frac{1}{2}}$-semistable for every $\alpha>0$ when $d=3,4$ and $\sigma_{\alpha, -\frac{1}{2}}$-semistable for $\alpha>\frac{1}{\sqrt{20}}$ when $d=5$. Hence $A$ is a $\mu$-semistable sheaf by \cite[Lemma 2.7]{BMS}. By \cite{Li15}, after taking $\alpha \to 0$, \cite[Conjecture 4.1]{BMS} holds for $A$ and $\alpha=0, \beta=-\frac{1}{2}$. And when $d=5$ we take $\alpha=\frac{1}{\sqrt{20}}, \beta=-\frac{1}{2}$.
	Thus we have $m\leq 1$ when $d=3,4$, and $m=0$ when $d=5$. Assume $m=1$, then we apply $\mathrm{Hom}(\mathcal{O}_Y,-)$ to $A[1]\to F\to B$ and obtain $\mathrm{hom}(\mathcal{O}_Y,A[2])=1$ and $\mathrm{hom}(\mathcal{O}_Y,A[i])=0$ for $i\neq 2$. By Serre duality we have $\mathrm{hom}(A, \mathcal{O}(-2)[1])=1$. But using Lemma \ref{no wall intersect -1}, a similar slope-comparison argument in Lemma \ref{inst-in-ku} shows that $\mathrm{hom}(A, \mathcal{O}(-2)[1])=0$, which makes a contradiction.

    Thus we conclude that $F=A[1]$, hence $F[-1]$ is $\sigma_{\alpha, -\frac{1}{2}}$-semistable  for every $\alpha >0$. By Lemma \ref{mu-ss-g-ss}, $F[-1]=A$ is actually a Gieseker-semistable sheaf. If $F$ is $\sigma(\alpha, -\frac{1}{2})$-stable, it can not be an extension of ideal sheaves of lines up to some shifts. Hence $F[-1]$ is Gieseker-stable in this case. 
    
    When $d=3,4$, by Lemma \ref{ins-prop} we know that  $F[-1]$ is an instanton sheaf. When $d=5$, from Lemma \ref{no wall intersect -1} we know that the only wall intersect with $\beta=-\frac{1}{2}$ is realized by $\mathcal{O}(-1)[2]$, and this is not an actual wall for $F$ by assumption. Since $\mu_{\alpha, -\frac{1}{2}}(\mathcal{O}(-1)[1])=\alpha^2-\frac{1}{4}$ and $\mu_{\alpha, -\frac{1}{2}}(F[-1])=-\frac{3}{20}-\alpha^2$, if we take $\alpha<\frac{1}{\sqrt{20}}$, we obtain $F, \mathcal{O}(-1)[2]\in \mathrm{Coh}^0_{\alpha, -\frac{1}{2}}(Y)$  and $\mu^0_{\alpha, -\frac{1}{2}}(F)>\mu^0_{\alpha, -\frac{1}{2}}(\mathcal{O}(-1)[2])$, which implies $\mathrm{Hom}(F[-1], \mathcal{O}(-1)[1])=0$. Hence $F[-1]$ is an instanton sheaf.
    
    Now assume that $\mathcal{C}$ is an actual wall. Then $F$ becomes unstable when $\beta_0\to -\frac{1}{2}$, and $F$ is strictly $\sigma^0_{\alpha_0, \beta_0}$-semistable and there's a sequence in
    \[0\to P\to F\to Q\to 0\]
    $\mathrm{Coh}^0_{\alpha_0, \beta_0}(Y)$, where $P,Q$ are $\sigma^0_{\alpha_0, \beta_0}$-semistable with the same slope $+\infty$. From Lemma \ref{infty-wall} we know  that $P\cong \mathcal{O}(-1)^{\oplus a}[2]$.
    
    Therefore, we have $\chi(P,Q)\neq 0$ and  $\mathrm{Hom}(P,Q[1])\neq 0$, since $P,Q$ are in the same heart and by Serre duality. So we can define an object $G$ as the extension
	 \[0\to Q\to G\to P\to 0\]
	 in $\mathrm{Coh}^0_{\alpha_0, \beta_0}(Y)$. Now $G$ is $\sigma^0_{\alpha, \beta}$-semistable when $\beta\to -\frac{1}{2}$. The argument in previous cases shows that $G[-1]$ is an instanton sheaf. Thus we have $\mathrm{Hom}(G[-1], \mathcal{O}(-1)[1])=\mathrm{Hom}(G, \mathcal{O}(-1)[2])=0$ and makes a contradiction.
	 
\end{proof}

Now we construct isomorphisms between moduli spaces:

\begin{theorem} \label{ins-iso}
Let $Y=Y_3,Y_4,Y_5$.	The projection functor $\mathrm{pr}$ induces an  isomorphism
	\[M^{inst}_Y \xrightarrow{\cong} \mathcal{M}^{ss}_{\sigma} (\mathcal{K}u(Y), 2L-2)\]
	for every Serre-invariant stability condition $\sigma$ on $\mathcal{K}u(Y)$.
\end{theorem}

\begin{proof}
	
	
	By Proposition \ref{all-serre-invariant}, we can assume $\sigma=\sigma(\alpha, -\frac{1}{2})$ for $0<\alpha<\frac{1}{2}$. Using Proposition \ref{sheaf-stable} and the GIT construction of $M^{ss}_Y(2,0,2)$, a similar argument in \cite[Section 5]{qin2019moduli} and \cite[Section 5]{qinv5} shows that the projection functor $\mathrm{pr}$ induces a morphism $M^{inst}_Y \to \mathcal{M}^{ss}_{\sigma} (\mathcal{K}u(Y), 2L-2)$. 
	
	By Proposition \ref{sheaf-stable} and Proposition \ref{stable-implies-sheaf}, we know that this morphism is bijective on closed points. Since $\mathrm{pr}$ effects nothing on $E$, we know that this morphism is etale.  Thus the projection functor induces a bijective etale morphism, which is an isomorphism.

\end{proof}



\begin{theorem} \label{ins-embd}
Let $Y=Y_1,Y_2$.	The projection functor $\mathrm{pr}$ induces an open immersion
\[M^{inst}_Y \hookrightarrow \mathcal{M}^{ss}_{\sigma} (\mathcal{K}u(Y), 2L-2)\]
for every Serre-invariant stability condition $\sigma\in \mathcal{K}$.
\end{theorem}

\begin{proof}
By \cite[Proposition 3.6]{PY20}, we can assume $\sigma=\sigma(\alpha, -\frac{1}{2})$. Using Proposition \ref{sheaf-stable} and the GIT construction of $M^{ss}_Y(2,0,2)$, a similar argument in \cite[Section 5]{qin2019moduli} and \cite[Section 5]{qinv5} shows that the projection functor $\mathrm{pr}$ induces a morphism $M^{inst}_Y \to \mathcal{M}^{ss}_{\sigma} (\mathcal{K}u(Y), 2L-2)$.

By Proposition \ref{sheaf-stable}, we know that this is injective on closed points. Since $\mathrm{pr}$ effects nothing on $E$, we know that this morphism is etale. Thus the projection functor induces an injective etale morphism, which is an open immersion.
\end{proof}

\subsection{All Serre-invariant stability conditions}

In this subsection, we assume $Y:=Y_d$ for $d=3,4,5$. We can consider all Serre-invariant stability conditions, not only for $\sigma\in \mathcal{K}$.

\begin{lemma} \label{1-5-heart}
	Let $F\in \mathcal{K}u(Y)$ be an object with $[F]=2v=2-2L$ and $\mathrm{hom}(F,F)=1, \mathrm{ext}^1(F,F)=5, \mathrm{ext}^2(F,F)=\mathrm{ext}^3(F,F)=0$. Then $F$ is in the heart of every Serre-invariant stability on $\mathcal{K}u(Y)$ up to some shifts.
\end{lemma}

\begin{proof}
	As in \cite[Lemma 4.5]{BMMS}, we consider the spectral sequence for objects in $\mathcal{K}u(Y)$ whose second page is given by 
	\[E^{p,q}_2=\bigoplus_i \mathrm{Hom}^p(\mathcal{H}^i(F), \mathcal{H}^{i+q}(F))\Rightarrow \mathrm{Hom}^{p+q}(F,F)\]
	where the cohomology is taken with respect to the heart.
	
	First we assume $d=3$. Then $E^{1,q}_2=E^{1,q}_{\infty}$ by Lemma \ref{heart-dim}. Thus $\mathrm{ext}^1(*,*)\geq 2$ for $*\in \mathcal{K}u(Y)$. If we take $q=0$, we obtain
	\[5=\mathrm{ext}^1(F,F)\geq \sum_i \mathrm{ext}^1(\mathcal{H}^i(F),\mathcal{H}^i(F))\geq 2r\]
	where $r$ is the number of non-zero cohomology objects of $F$.
	
	If $r=1$. then $F$ is already in the heart up to some shifts. Otherwise if $r=2$, we denote these two cohomology objects by $M, N$. Thus we have $\mathrm{ext}^1(M,M)=\mathrm{ext}^1(N,N)=2$ or $\mathrm{ext}^1(M,M)=2, \mathrm{ext}^1(N,N)=3$ or $\mathrm{ext}^1(M,M)=3, \mathrm{ext}^1(N,N)=2$. But $\chi(M,M)\leq -1$ and $\chi(N,N)\leq -1$, hence if $\mathrm{ext}^1(M,M)=2, \mathrm{ext}^1(N,N)=3$, then $-1\leq \chi(M,M)\leq -1$ and $-2\leq \chi(N,N)\leq -1$. Since there's no (-2)-classes in $\mathcal{K}u(Y_3)$, we have $\chi(M,M)=\chi(N,N)=-1$. In other two cases we also have $\chi(M,M)=\chi(N,N)=-1$, i.e. $M$ and $N$ are always (-1)-classes. 
	Now since $\mathrm{ch}(F)=2-2L=2v$, by classifications of $(-1)$-classes we know $\mathrm{ch}(M)=\mathrm{ch}(N)=1-L=v$.
	Thus the second page is:
	
	$$E^{p,q}_2= \begin{array}{|ccc}
	     0 & 0 & 0\\
	     \mathrm{Hom}(N,M) & \mathrm{Ext}^1(N,M) & \mathrm{Ext}^2(N,M) \\
	     \mathrm{Hom}(M,M)\oplus \mathrm{Hom}(N,N) & \mathrm{Ext}^1(M,M)\oplus \mathrm{Ext}^1(N,N) & 0 \\ \hline
	     \mathrm{Hom}(M,N) &  \mathrm{Ext}^1(M,N) & \mathrm{Ext}^2(M,N)
	\end{array}$$

The table of dimension of page 2 is:

$$\dim E^{p,q}_2= \begin{array}{|ccc}
	     0 & 0 & 0\\
	     a & b & c \\
	     2 & 4 & 0 \\ \hline
	     d &  e & f
	\end{array}$$

From the convergence of spectral sequence, we know $c=0$, $b=0$. And from $\chi(N,M)=-1$ we know $a-b+c=-1$, which implies $a=-1$, thus makes a contradiction.

Next we assume $d=4,5$. In both cases the heart of a stability condition on $\mathcal{K}u(Y)$ has homological dimension 1, hence $E^{p,q}_2=E^{p,q}_{\infty}$. If we take $q=0$, we obtain
	\[5=\mathrm{ext}^1(F,F)\geq \sum_i \mathrm{ext}^1(\mathcal{H}^i(F),\mathcal{H}^i(F))\geq r\]
	If $r=1$ then we done. Otherwise if $r\geq 2$, then this is also impossible since $1=\mathrm{hom}(F,F)=\dim E^{0,0}_2+\dim E^{1,-1}_2\geq \dim E^{0,0}_2= 2$ from the convergence of spectral sequence, which makes a contradiction.
\end{proof}

From some elementary computations, we have:

\begin{lemma} \label{phi-equal}
 Let $A,B\in \mathcal{K}u(Y)$. Assume $\mathrm{ch}(A)=a_1v+b_1w$ and $\mathrm{ch}(B)=a_2v+b_2w$ for $a_i, b_i\in \mathbb{Z}$, then $\mu_{\alpha, -\frac{1}{2}}(A)=\mu_{\alpha, -\frac{1}{2}}(B)$ if and only if in the following cases:
 \begin{enumerate}
    \item $a_1=a_2=0$ or $b_1=b_2=0$
    
    \item $a_1, a_2, b_1, b_2\neq 0$, and $\frac{a_1}{b_1}=\frac{a_2}{b_2}$.
\end{enumerate}
\end{lemma}


\begin{proposition} \label{1-5-serre-stable}
	Let $F\in \mathcal{K}u(Y)$ be an object with $\mathrm{ch}(F)=2-2L$ and $\mathrm{hom}(F,F)=1, \mathrm{ext}^1(F,F)=5, \mathrm{ext}^2(F,F)=\mathrm{ext}^3(F,F)=0$, then $F$ is $\sigma$-stable for every Serre-invariant satbility condition on $\mathcal{K}u(Y)$. In particular, every $E\in M_Y(2,0,2)$ is $\sigma$-stable.
	
\end{proposition}

\begin{proof}
	By Lemma \ref{1-5-heart}, without loss of generality we can assume that $F\in \mathcal{A}_{\sigma}$. 
	
If $F$ is not $\sigma$-semistable, then there's an exact triangle
$$A\to F\to B$$
such that $A,B$ are both $\sigma$-semistable with $\phi(A)>\phi(B)$. Hence we have  $\mathrm{Hom}(A,B)=0$ and  $\mathrm{ext}^2(B,A)=\mathrm{ext}^2(A,A)=\mathrm{ext}^2(B,B)=0$. 

First we assume $d=3$. By Lemma \ref{mk-lem} we have
\[\mathrm{ext}^1(A,A)+\mathrm{ext}^1(B,B)\leq \mathrm{ext}^1(F,F)=5\]
By Lemma \ref{py-ext-2}, the only possible cases are $\mathrm{ext}^1(A,A)=\mathrm{ext}^1(B,B)=2$ or $\mathrm{ext}^1(A,A)=2, \mathrm{ext}^1(B,B)=3$ or $\mathrm{ext}^1(A,A)=3, \mathrm{ext}^1(B,B)=2$. But $\chi(A,A)\leq -1$ and $\chi(B,B)\leq -1$, hence if $\mathrm{ext}^1(A,A)=2, \mathrm{ext}^1(B,B)=3$, then $-1\leq \chi(A,A)\leq -1$ and $-2\leq \chi(B,B)\leq -1$. Since there's no (-2)-classes in $\mathcal{K}u(Y_3)$, we have $\chi(A,A)=\chi(B,B)=-1$. In other two cases we also have $\chi(A,A)=\chi(B,B)=-1$, i.e. $A$ and $B$ are always (-1)-classes. Now since $\mathrm{ch}(F)=2-2L=2v$, by classfication of $(-1)$-classes we know $\mathrm{ch}(A)=\mathrm{ch}(B)=1-L=v$. But in this case we have $\phi(A)=\phi(B)$, which is a contradiction. 

When $d=4,5$, the heart $\mathcal{A}_{\sigma}$ of  $\mathcal{K}u(Y)$ has homological dimension 1, then it's easy to make contradictions by spectral sequence in \cite[Proposition 4.16]{zhang2020bridgeland} and Mukai Lemma. First we assume $d=5$, then the only possible cases are $\mathrm{ext}^1(A,A)=\mathrm{ext}^1(B,B)=2$ or $\mathrm{ext}^1(A,A)=2, \mathrm{ext}^1(B,B)=3$ or $\mathrm{ext}^1(A,A)=3, \mathrm{ext}^1(B,B)=2$.
	
	From the triangle $B[-1]\to A\to F$ we have a spectral sequence which degenerates at $E_3$ converging to $\mathrm{Ext}^*(F,F)$ with $E_1$-page being
	\[E^{p,q}_1= \left\{
	\begin{aligned}
	\mathrm{Ext}^q(A,B[-1])=\mathrm{Ext}^{q-1}(A,B) & , ~ p=-1 \\
	\mathrm{Ext}^q(A,A)\oplus \mathrm{Ext}^q(B,B) & , ~ p=0 \\
	\mathrm{Ext}^q(B[-1], A)=\mathrm{Ext}^{q+1}(B,A) & , ~ p=1 \\
	0 &, ~  p\notin [-1,1]
	\end{aligned}
	\right.\]
	
	When $\mathrm{ext}^1(A,A)+\mathrm{ext}^1(B,B)=4$, 
    the dimension of $E_1$ page looks like:
	
	$$\dim E^{p,q}_1= \begin{array}{cc|cc}
	0 & 0 & 0 & 0\\
	0 & 0 & 0 & 0\\
	0 & c & 0 & 0\\
	0 & 0 & 4 & 0 \\
	0 & 0 & 2 & b \\ \hline
	0 & 0 &  0 & a
	\end{array}$$
	
	The dimension of $E_2$ page looks like:
	
	$$\dim E^{p,q}_2= \begin{array}{cc|cc}
	0 & 0 & 0 & 0\\
	0 & 0 & 0 & 0\\
	0 & c & 0 & 0\\
	0 & 0 & 4 & 0 \\
	0 & 0 & x & y \\ \hline
	0 & 0 &  0 & a
	\end{array}$$
	Hence we have $x+a=1, y+c=1, x-y=2-b$ and all variables are non-negative. Whenever $x=0$ or $x=1$, we have $\chi(A,B)=-1, \chi(B,A)=-1$ or $\chi(A,B)=0, \chi(B,A)=-2$. Assume $\mathrm{ch}(A)=a_1v+b_1w$ and $\mathrm{ch}(B)=a_2v+b_2w$, we have $a_1+a_2=2$, $b_1+b_2=0$ and $a_i, b_i\in \mathbb{Z}$. When $\chi(A,B)=0, \chi(B,A)=-2$, we have $-a_2(a_1+b_1)+b_2((1-d)a_1-db_1)=0$ and $-a_1(a_2+b_2)+b_1((1-d)a_2-db_2)=-2$, which has no solution under our assumptions. When $\chi(A,B)=-1, \chi(B,A)=-1$, we have $-a_2(a_1+b_1)+b_2((1-d)a_1-db_1)=-1$ and $-a_1(a_2+b_2)+b_1((1-d)a_2-db_2)=-1$, which has the only solution $a_1=a_2=1, b_1=b_2=0$.

	When $\mathrm{ext}^1(A,A)+\mathrm{ext}^1(B,B)=5$, since there's no $(-2)$-class, 
	the dimension of $E_1$ page looks like:
	
	$$\dim E^{p,q}_1= \begin{array}{cc|cc}
	0 & 0 & 0 & 0\\
	0 & 0 & 0 & 0\\
	0 & c & 0 & 0\\
	0 & 0 & 5 & 0 \\
	0 & 0 & 3 & b \\ \hline
	0 & 0 &  0 & a
	\end{array}$$
	
	The dimension of $E_2$ page looks like:
	
	$$\dim E^{p,q}_2= \begin{array}{cc|cc}
	0 & 0 & 0 & 0\\
	0 & 0 & 0 & 0\\
	0 & c & 0 & 0\\
	0 & 0 & 5 & 0 \\
	0 & 0 & x & y \\ \hline
	0 & 0 &  0 & a
	\end{array}$$
	Hence we have $c=y=0, x+a=1, x=3-b$ and all variables are non-negative. Whenever $x=0$ or $x=1$, we have $\chi(A,B)=0, \chi(B,A)=-2$. Assume $\mathrm{ch}(A)=a_1v+b_1w$ and $\mathrm{ch}(B)=a_2v+b_2w$, we have $a_1+a_2=2$, $b_1+b_2=0$ and $a_i, b_i\in \mathbb{Z}$. Since $\chi(A,B)=0, \chi(B,A)=-2$, we have $-a_2(a_1+b_1)+b_2((1-d)a_1-db_1)=0$ and $-a_1(a_2+b_2)+b_1((1-d)a_2-db_2)=-2$, which has no solution under our assumptions.

	Next we assume $d=4$.
	In this case we must have $\chi(A,A)+\chi(B,B)=-4$, hence $\mathrm{ext}^1(A,A)+\mathrm{ext}^1(B,B)=4+\mathrm{hom}^1(A,A)+\mathrm{hom}^1(B,B)\geq 6$, which contradicts with Mukai lemma.

Now we are going to show that  $F$ is $\sigma$-stable. By \cite[Proposition 3.6]{PY20} we can assume $\sigma=\sigma(\alpha, -\frac{1}{2})$.  Assume $F$ is strictly $\sigma$-semistable. If $F$ has more than one non-isomorphic stable factors, then there's an exact sequence
$$0\to A\to F\to B\to 0$$
where $\phi(A)=\phi(B)=\phi(F)$ and $A,B$ are $\sigma$-semistable with $\mathrm{Hom}(A,B)$=0. When $d=3,5$, by Mukai Lemma we know $\chi(A,A)=\chi(B,B)=-1$. Since $\mathrm{ch}(A)+\mathrm{ch}(B)=\mathrm{ch}(F)=2v$, we know that only possible case is $\mathrm{ch}(A)=\mathrm{ch}(B)=v$. When $d=4$, since there's only one orbit in $\mathrm{Stab}(\mathcal{K}u(Y_4))$, we can assume $\sigma=\sigma(\alpha, -\frac{1}{2})$. Then by Lemma \ref{phi-equal} we know that only possible case is $\mathrm{ch}(A)=\mathrm{ch}(B)=v$.
But as show in \cite[Proposition 4.6]{PY20}, this implies that $A\cong I_L[2k]$ and $B\cong I_{L'}[2k]$ for some lines $L,L'$ on $Y$ and $k\in \mathbb{Z}$. If all stable factors of $F$ are isomorphic to $S$, then $\chi(F,F)=n^2\chi(S,S)=-4$, which implies $n=2$. This means $F$ is also an extension of $S\cong I_L[2k]$. But both cases contradict with our assumptions by \cite[Lemma 4.3]{D}.

\end{proof}



By \cite[Remark 5.14]{PY20}, for every pair of Serre-invariant stability conditions $\sigma_1=(\mathcal{P}_1, Z_1)$ and $\sigma_2=(\mathcal{P}_2, Z_2)$ on $\mathcal{K}u(Y)$ we have $Z_1=T\circ Z_2$ for some $T=(t_{i j})_{1\leq i,j\leq 2}\in \mathrm{GL}^+_2(\mathbb{R})$. Thus an elementary computation shows that:

\begin{lemma} \label{same-slope}
Let $\sigma_1=(\mathcal{P}_1, Z_1)$ and $\sigma_2=(\mathcal{P}_2, Z_2)$ be two Serre-invariant stability conditions on $\mathcal{K}u(Y)$. Let $E,E'\in \mathcal{K}u(Y)$ be any two objects, then $\mu_1(E)=\mu_1(E')$ if and only if $\mu_2(E)=\mu_2(E')$; And $\mu_1(E)>\mu_1(E')$ if and only if $\mu_2(E)>\mu_2(E')$.
\end{lemma}

	

\begin{lemma} \label{non-negative}
Let $\sigma_1=(\mathcal{P}_1, Z_1)$ be a Serre-invariant stability condition on $\mathcal{K}u(Y)$ and $E\in \mathcal{P}_1(\phi_1(E))$ be a $\sigma_1$-stable object with $[E]=av$ for $a\in \mathbb{Z}$. Then for any $\sigma_1$-stable object $F\in \mathcal{P}_1(\phi_1(E))$ we have $[F]=bv$ for $b\in \mathbb{Z}$ and $ab>0$.
\end{lemma}

\begin{proof}
       When $F\cong E$, this is clear. Assume $E$ and $F$ are not isomorphic. Then by Lemma \ref{same-slope} and Lemma \ref{phi-equal} we know $[F]=bv$ for $b\in \mathbb{Z}$. It's clear that $ab\neq 0$. If $ab<0$, then $\chi(F,E)=\mathrm{hom}(F,E)-\mathrm{ext}^1(F,E)=-ab>0$. This means $\mathrm{hom}(F,E)>0$. But this is impossible since $E$ and $F$ are both stable and not isomorphic.
\end{proof}

\begin{proposition} \label{all-serre-invariant}
Let $\sigma_1=(\mathcal{P}_1, Z_1)$ and $\sigma_2=(\mathcal{P}_2, Z_2)=\sigma(\alpha, -\frac{1}{2})$ be two Serre-invariant stability conditions on $\mathcal{K}u(Y)$. Let $E\in \mathcal{K}u(Y)$ be an object with $[E]=\pm2v$. Then $E\in \mathcal{P}_1((0,1])$ is $\sigma_1$-(semi)stable if and only if $E[m]\in \mathcal{P}_2((0,1])$ is $\sigma_2$-(semi)stable for $m=m(\sigma_1)\in \mathbb{Z}$.
\end{proposition}

\begin{proof}
Without loss of generality, we can assume $[E]=-2v$.
As shown in \cite{PY20}, if $F\in \mathcal{K}u(Y)$ is a $\sigma_1$-stable object with $[F]=-v$, then $F$ is isomorphic to an ideal sheaf of line up to some shifts, and this shifts only depend on $\sigma_1$; And if $F\in \mathcal{K}u(Y)$ is a $\sigma_2$-stable object with $[F]=-v$, then $F\cong I_L[1]$ where $L$ is a line on $Y$. We can assume $I_L[1-m]\in \mathcal{P}_1((0,1])$, where $m=m(\sigma_1)\in \mathbb{Z}$ only depends on $\sigma_1$. Thus if $F\in \mathcal{K}u(Y)$ is an object with $[F]=-v$, then $F\in \mathcal{P}_1((0,1])$ is $\sigma_1$-stable if and only if $F[m]\in \mathcal{P}_2((0,1])$ is $\sigma_2$-stable.



 By Proposition \ref{1-5-serre-stable}, $E\in \mathcal{P}_1((0,1])$ is $\sigma_1$-stable if and only if $E$ is  $\sigma_2$-stable up to some shifts.
To show $E[m]$ is in the heart of $\sigma_2$, we can do as in  \cite[Section 5.1]{APR}: Let $F=I_L[1-m]$. Since $\mathrm{Hom}(F,E[1])=\mathrm{Ext}^1(F,E)\neq 0$, we have $\phi_2(F)<\phi_2(E)+1$. Since $\mathrm{Hom}(E,F[1])=\mathrm{Ext}^1(E,F)\neq 0$, we have $\phi_2(E)<\phi_2(F)+1$. Thus $-1<\phi_2(E)-\phi_2(F)<1$. But by Lemma \ref{same-slope} we know $\mu_2(E)=\mu_2(F)$, thus $\phi_2(E)-\phi_2(F)\in \mathbb{Z}$, which shows $\phi_2(E)=\phi_2(F)$. Thus since $F[m]\in \mathcal{P}_2((0,1])$, we have $E[m]\in \mathcal{P}_2((0,1])$.

Now we assume $E\in \mathcal{P}_1((0,1])$ is strictly $\sigma_1$-semistable, let $E_1,...,E_n$ be the JH-factors under $\sigma_1$ of $E$. Then $E_i$ are $\sigma_1$-stable with $\phi_1(E_1)=...=\phi_1(E_n)$. Since $\sum \mathrm{ch}(E_i)=\mathrm{ch}(E)=-2v$, by Lemma \ref{phi-equal}, Lemma \ref{same-slope} and Lemma \ref{non-negative} we have $n=2$ and $\mathrm{ch}(E_1)=\mathrm{ch}(E_2)=-v$. This means $E_1\cong I_{L_1}[1-m]$ and $E_2\cong I_{L_2}[1-m]$, where $L_1,L_2$ are two lines on $Y$. Thus we know that $E[m]$ is an extension of $I_{L_1}[1]$ and $I_{L_2}[1]$, hence is $\sigma_2$-semistable.

Conversely, if $E$ is $\sigma_2$-stable, then by Proposition \ref{1-5-serre-stable} and the same argument above shows that $E[-m]$ is $\sigma_1$-stable. If $E$ is strictly $\sigma_2$-semistable, then $E$ is an extension of two ideal sheaves of lines up to shift by 1. Thus we also have $E[-m]$ is $\sigma_1$-semistable.
\end{proof}

\section{Isomorphisms between Bridgeland moduli spaces of $(-4)$-classes} \label{sec-8}

Let $X=X_{4d+2}$ for $d=3,4,5$. In this section we are going to show that $\mathcal{M}^{ss}_{\sigma}(\mathcal{K}u(Y), 2-2L)\cong \mathcal{M}^{ss}_{\sigma'}(\mathcal{A}_X, 2-4L)$ and $\mathcal{M}_{\sigma}(\mathcal{K}u(Y), 2-2L)\cong \mathcal{M}_{\sigma'}(\mathcal{A}_X, 2-4L)$.


Almost the same as Lemma \ref{1-5-heart} and Proposition \ref{1-5-serre-stable}, but using Corollary \ref{406} instead of \cite[Proposition 4.6]{PY20}, we have:

\begin{lemma} \label{1-5-heart-X}
	Let $F\in \mathcal{A}_X$ be an object with $[F]=2s=2-4L$ and $\mathrm{hom}(F,F)=1, \mathrm{ext}^1(F,F)=5, \mathrm{ext}^2(F,F)=\mathrm{ext}^3(F,F)=0$. Then up to some shifts $F$ is in the heart of every Serre-invariant stability on $\mathcal{A}_X$.
\end{lemma}




\begin{proposition} \label{1-5-serre-stable-X}
	Let $F\in \mathcal{A}_X$ be an object with $\mathrm{ch}(F)=2-4L$ and $\mathrm{hom}(F,F)=1, \mathrm{ext}^1(F,F)=5, \mathrm{ext}^2(F,F)=\mathrm{ext}^3(F,F)=0$, then $F$ is $\sigma'$-stable with respect to every Serre-invariant satbility condition $\sigma'$ on $\mathcal{A}_X$.
	
\end{proposition}

Using the equivalence in \cite{KPS}, Proposition \ref{all-serre-invariant} holds for $\mathcal{A}_{X_{4d+2}}$ and $d=3,4,5$ if we replace $\sigma(\alpha, -\frac{1}{2})$ by $\Phi(\sigma(\alpha, -\frac{1}{2}))$. In summary, we have:

\begin{corollary} \label{bi-ins}
Let $(Y,X)\in Z_d\subset \mathcal{MF}^2_{d}\times \mathcal{MF}^1_{4d+2}$ for $d=3,4,5$. Let $E$ be an object in $\mathcal{K}u(Y)$ with $[E]=2-2L$. Let $\sigma$ and $\sigma'$ be two Serre-invariant stability conditions on $\mathcal{K}u(Y)$ and $\mathcal{A}_X$ respectively. Then $E\in \mathcal{P}_{\sigma}((0,1])$ is $\sigma$-(semi)stable if and only if $\Phi(E)[m]\in \mathcal{P}_{\sigma'}((0,1])$ is $\sigma'$-(semi)stable for $m=m(\sigma, \sigma')\in \mathbb{Z}$.
\end{corollary}

\begin{theorem} \label{iso-4}
Let $(Y,X)\in Z_d\subset \mathcal{MF}^2_{d}\times \mathcal{MF}^1_{4d+2}$ for $d=3,4,5$.	Then the equivalence $\Phi: \mathcal{K}u(Y)\cong \mathcal{A}_X$ induces an isomorphism between moduli spaces:
\[s: \mathcal{M}^{ss}_{\sigma}(\mathcal{K}u(Y), 2-2L) \xrightarrow{\cong} \mathcal{M}^{ss}_{\sigma'}(\mathcal{A}_X, 2-4L)\]
such that the restriction
	\[\mathcal{M}_{\sigma}(\mathcal{K}u(Y), 2-2L) \xrightarrow{\cong} \mathcal{M}_{\sigma'}(\mathcal{A}_X, 2-4L)\]
	is also an isomorphism. Here $\sigma$ and $\sigma'$ are Serre-invariant stability conditions on $\mathcal{K}u(Y)$ and $\mathcal{A}_X$ respectively.
\end{theorem}

\begin{proof}
	From Theorem \ref{ins-iso} we know $M^{inst}_Y \cong \mathcal{M}^{ss}_{\sigma} (\mathcal{K}u(Y), 2L-2)$. Since the equivalences in \cite{KPS} map the numerical class $2(1-L)$ to $2(1-2L)$, hence using Corollary \ref{bi-ins}, a similar argument in Theorem \ref{ins-iso} shows that $\Phi$ induces a morphism $s$. This is a proper morphism since moduli spaces on both sides are proper.
	
	From Corollary \ref{bi-ins} we know that $s$ is bijective at closed points and $\mathcal{M}_{\sigma}(\mathcal{K}u(Y), 2-2L)=s^{-1}(\mathcal{M}_{\sigma'}(\mathcal{A}_X, 2-4L))$. Since $s$ is induced by equivalence, it is etale. Now we know $s$ is a bijective etale proper morphism, which is an isomorphism.
\end{proof}




\begin{corollary} \label{all moduli spaces of (-4)-class isom}
Let $b$ be any $(-4)$-class in $\mathcal{N}(\mathcal{K}u(Y_d))$ and $b'$ be any $(-4)$-class in $\mathcal{N}(\mathcal{A}_{X_{4d+2}})$. Then the Bridgeland moduli spaces $\mathcal{M}^{ss}_{\sigma}(\mathcal{K}u(Y_d),b)\cong\mathcal{M}^{ss}_{\sigma'}(\mathcal{A}_{X_{4d+2}},b')$ for any Serre invariant stability condition $\sigma,\sigma'$ and for each $d=3,5$. 
\end{corollary}

\begin{proof}
       When $d=3,5$, note that $(x_0,y_0)$ is an integer solution of $x^2+dxy+dy^2=4$ if and only if $(x_0/2,y_0/2)$ is an integer solution of $x^2+dxy+dy^2=1$, Thus the result follows from a similar argument in Corollary \ref{all moduli spaces of (-1)-class isom}.
\end{proof}

\section{Moduli space $M^{ss}_X(2,0,4)$ on $X_{4d+2}$ as Bridgeland moduli space} \label{sec-9}

In this section, we are going to show that there's an isomorphism between $M^{ss}_X(2,0,4)$ and $\mathcal{M}^{ss}_{\sigma'}(\mathcal{A}_X, 2-4L)$
for $X=X_{14}, X_{18}$ and $X_{22}$.

As shown in \cite{BCFacm}, there're following classes of sheaves in $M^{ss}_X(2,0,4)$:
	
	\begin{enumerate}
		\item $E$ is strictly Gieseker-semistable if and only if $E$ is an extension of ideal sheaves of conics.
		
		\item $E$ is Gieseker-stable and locally free.
		
		\item $E$ is Gieseker-stable but not locally free, and fit into an exact sequence
		\[0\to I_C\to E\to I_L\to 0\]
		where $C$ is a cubic and $L$ is a line. In this case $\chi(I_C)=0$, $\mathrm{ch}(I_C)=1-3L+\frac{1}{2}P$.
		
		\item $E$ is Gieseker-stable but not locally free, and fit into an exact sequence
		\[0\to I_C\to E\to I_x\to 0\]
		where $C$ is a quartic and $x\in C$ is a point. In this case $\chi(I_C)=0$, $\mathrm{ch}(I_C)=1-4L+P$.
		
		\item $E$ is Gieseker-stable but not locally free, and fit into an exact sequence
		\[0\to I_C\to E\to \mathcal{O}_X\to 0\]
		where $C$ is of degree $4$ with $\chi(\mathcal{O}_C)=2$. In this case $\chi(I_C)=-1$, $\mathrm{ch}(I_C)=1-4L$.
	\end{enumerate}
	
	When $E$ is strictly Gieseker-semistable, it's clear from classifications above that $E\in \mathcal{A}_X$ and is semistable with respect to every Serre-invariant stability condition. When $E$ is Gieseker-stable, this is also true.
	
	\begin{lemma} \label{204-in-ku}
		Let $E\in M_X(2,0,4)$, then $h^i(E)=0$ and $\mathrm{ext}^i(\mathcal{E}^{\vee}, E)=0$ for every $i$. Thus $E\in \mathcal{A}_X$.
	\end{lemma}

\begin{proof}
	By \cite[Lemma 4.3]{BCFacm} we know $h^2(E)=0$. From the stability and Serre duality, we have   $h^3(E)=\mathrm{ext}^3(\mathcal{O}_X, E)=\mathrm{hom}(E, \mathcal{O}_X(-1))=0$ and $h^0(E)=\mathrm{hom}(\mathcal{O}_X, E)=0$. Since $\chi(E)=0$, we know $h^1(E)=0$. Thus $h^i(E)=0$ for every $i$.
	
	To prove $\mathrm{ext}^i(\mathcal{E}^{\vee}, E)=0$, we assume $d=3$ for simplicity, but other cases are almost the same. 
	
	First we assume that $E$ is in case (3). From stability we have $\mathrm{hom}(\mathcal{E}^{\vee}, I_C)=\mathrm{ext}^3(\mathcal{E}^{\vee}, I_C)=\mathrm{hom}(\mathcal{E}^{\vee}, I_L)=\mathrm{ext}^3(\mathcal{E}^{\vee}, I_L)=0$. From $\mathcal{E}|_L=\mathcal{O}_L\oplus \mathcal{O}_L(-1)$ we know  $\mathrm{ext}^1(\mathcal{E}^{\vee}, I_L)=1$ and $\mathrm{ext}^2(\mathcal{E}^{\vee}, I_L)=0$. Since $\chi(\mathcal{E}^{\vee}, I_C)=1$, we only need to show that $\mathrm{ext}^1(\mathcal{E}^{\vee}, I_C)=0$, that is $h^0(\mathcal{E}|_C)=0$. 
	
	If $E$ is in case (4), from stability we have $\mathrm{hom}(\mathcal{E}^{\vee}, I_C)=\mathrm{ext}^3(\mathcal{E}^{\vee}, I_C)=\mathrm{hom}(\mathcal{E}^{\vee}, I_x)=\mathrm{ext}^3(\mathcal{E}^{\vee}, I_x)=0$. From the long exact sequence we have  $\mathrm{ext}^1(\mathcal{E}^{\vee}, I_x)=2$. Since $\chi(\mathcal{E}^{\vee}, I_C)=2$, we only need to show that $\mathrm{ext}^1(\mathcal{E}^{\vee}, I_C)=0$, that is $h^0(\mathcal{E}|_C)=0$.
	
	If $E$ is in case (5), from stability we have $\mathrm{ext}^*(\mathcal{E}^{\vee}, \mathcal{O}_X)=0$ and $\mathrm{hom}(\mathcal{E}^{\vee}, I_C)=\mathrm{ext}^3(\mathcal{E}^{\vee}, I_C)=0$. Since $\chi(\mathcal{E}^{\vee}, I_C)=0$, we only need to show that  $\mathrm{ext}^1(\mathcal{E}^{\vee}, I_C)=0$, that is $h^0(\mathcal{E}|_C)=0$. Now for all these three cases, a similar argument in \cite[Lemma B.3.3]{KPS} shows that these cohomological groups are vanishing as we want.

	When $E$ is locally free, from stability and Serre duality we have $\mathrm{hom}(\mathcal{E}^{\vee}, E)=0$ and $\mathrm{ext}^3(\mathcal{E}^{\vee}, E)=0$.
	Since $\chi(\mathcal{E}^{\vee}, E)=0$, we only need to show $\mathrm{ext}^1(\mathcal{E}^{\vee}, E)=\mathrm{hom}(\mathcal{E}^{\vee}, E[1])=0$. To this end, 
	by \cite[Lemma 2.7 (ii)]{BMS} we know that $\mathcal{E}^{\vee}$ and $E[1]$ are in $\mathrm{Coh}^{\frac{1}{4}}(X)$ and are $\sigma_{\alpha, \frac{1}{4}}$-stable for $\alpha \gg 0$. But $\mu_{\alpha, \frac{1}{4}}(\mathcal{E}^{\vee})> \mu_{\alpha, \frac{1}{4}}(E[1])$ when $\alpha\gg 0$, which gives $\mathrm{hom}(\mathcal{E}^{\vee}, E[1])=\mathrm{ext}^1(\mathcal{E}^{\vee}, E)=0$.
	
\end{proof}

	Now we are going to show that $E\in M_X(2,0,4)$ is stable with respect to every Serre-invariant stability condition $\sigma'$ on $\mathcal{A}_X$.
	
	\begin{lemma} \label{204-stable}
	Let $E\in M_X(2,0,4)$, then $\mathrm{hom}(E,E)=1,\mathrm{ext}^1(E,E)=5, \mathrm{ext}^2(E,E)=\mathrm{ext}^3(E,E)=0$.
	\end{lemma}
	
	\begin{proof}
	  It's clear $\chi(E, E)=-4$. And $\mathrm{ext}^3(E,E)=\mathrm{hom}(E,E(-1))=0$ from stability. Hence we only need to show $\mathrm{ext}^2(E,E)=0$. 
	  
	  By \cite[Section 6]{BLMS}, $\sigma(\alpha, \beta)=\sigma_{\alpha, \beta}^0|_{\mathcal{A}_X}$ is a stability condition on $\mathcal{A}_X$ for every  $\beta<0$ and $0<\alpha$ with $-\beta, \alpha$ both sufficiently small. When $d=4,5$, the heart $\mathcal{A}(\alpha, \beta)$ has homological dimension 1. When $d=3$, from \cite{Pertusi2021serreinv} we know that these stability conditions are Serre-invariant. Thus by Lemma \ref{py-phase}, to show $\mathrm{ext}^2(E,E)=0$, we only need to show $E\in \mathcal{A}_X$ is $\sigma(\alpha, \beta)$-semistable.
	  
	  To this end, from \cite[Proposition 4.8]{bayer2020desingularization} we know that $E$ is $\sigma_{\alpha,\beta}$-semistable for all $\beta<0$ and $\alpha\gg 0$. Thus by locally-finiteness of walls, to show $E$ is $\sigma_{\alpha, \beta}$-semistable for every $\alpha>0$ and $\beta<0$ with $-\beta$ sufficiently small, we only need to show that there's no semicircle wall tangent with line $\beta=0$. If there's a such semicircle wall $\mathcal{C}$, then it's given by a sequence
	  \[0\to A\to E\to B\to 0\]
	  in $\mathrm{Coh}^{\beta}(X)$ for $(\alpha, \beta)\in \mathcal{C}$, such that $\mu_{\alpha,\beta}(A)=\mu_{\alpha,\beta}(E)=\mu_{\alpha,\beta}(B)$. Assume $\mathrm{ch}_{\leq 2}(A)=(a,bH, \frac{c}{d}H^2)$ for $a,b,c\in \mathbb{Z}$. By continuity, we have $\mu_{0,0}(A)=\mu_{0,0}(E)=\mu_{0,0}(B)$. Thus we have $b=0$. But as in Proposition \ref{wall-cross-1}, $\mathrm{ch}^{\beta}_1(A)\geq 0$ and $\mathrm{ch}^{\beta}_1(B)\geq 0$ for $\beta<0$ implies $a=0,1,2$, and $\mu_{\alpha,\beta}(A)=\mu_{\alpha,\beta}(E)=\mu_{\alpha,\beta}(B)<+\infty$ implies $a\neq 0,2$. Thus the only possible case is $\mathrm{ch}_{\leq 2}(A)=\mathrm{ch}_{\leq 2}(B)=(1,0,-2L)$, which effect nothing on stability of $E$ and are not semicircle walls. Therefore, there's no wall for $E$ with respect to $\sigma_{\alpha, \beta}$ on $\beta=-\epsilon$ where $\epsilon>0$ sufficiently small. Thus $E$ is $\sigma_{\alpha, \beta}$-semistable for every $\alpha>0$ and $\beta<0$ with $-\beta$ sufficiently small. Hence $E[1]\in \mathrm{Coh}^0_{\alpha, \beta}(X)$ is $\sigma^0_{\alpha, \beta}$-semistable, and $E[1]\in \mathcal{A}(\alpha, \beta)$ is $\sigma(\alpha, \beta)$-semistable for such $(\alpha, \beta)$.
	\end{proof}

	

	\begin{proposition} \label{acm-s-ss}
	Let $E\in M_X^{ss}(2,0,4)$. Then $E\in \mathcal{A}_X$ and if $E$ is Gieseker-(semi)stable, then $E$ is $\sigma'$-(semi)stable.
	\end{proposition}
	
	\begin{proof}
	The statement follows from Lemma \ref{204-in-ku}, Lemma \ref{204-stable}, Lemma \ref{1-5-heart-X} and Corollary \ref{bi-ins}.
	\end{proof}

	\begin{theorem} \label{acm-iso}
		Let $X=X_{14}, X_{18}$ or $X_{22}$. Then the projection functor induces an isomorphism
		\[s': M^{ss}_X(2,0,4)\xrightarrow{\cong} \mathcal{M}^{ss}_{\sigma'}(\mathcal{A}_X, 2-4L)\] 
		and the restriction gives an isomorphism
		\[M_X(2,0,4)\xrightarrow{\cong} \mathcal{M}_{\sigma'}(\mathcal{A}_X, 2-4L)\]
	for every Serre-invariant stability condition $\sigma'$ on $\mathcal{A}_X$.
		
	\end{theorem}

	\begin{proof}
       A similar argument in Theorem \ref{ins-iso} shows that the projection functor $\mathrm{pr}$ induces a morphism $s': M^{ss}_X(2,0,4)\to \mathcal{M}^{ss}_{\sigma'}(\mathcal{A}_X, 2-4L)$. We know that $M^{ss}_X(2,0,4)$ is projective and  $\mathcal{M}^{ss}_{\sigma'}(\mathcal{A}_X, 2-4L)$ is proper, hence $s'$ is also projective.  By Proposition \ref{acm-s-ss}, we have  $s'^{-1}(\mathcal{M}_{\sigma'}(\mathcal{A}_X, 2-4L))=M_X(2,0,4)$ and $s'$ is injective.
	
		From the fact that $\mathrm{pr}$ effects nothing on $E\in M^{ss}_X(2,0,4)$, we know $s'$ is etale. Thus $s'$ is an embedding. Since $\mathcal{M}^{ss}_{\sigma'}(\mathcal{A}_X, 2-4L)$ is irreducible and smooth by Lemma \ref{ins-prop}, Theorem \ref{ins-iso} and Theorem \ref{iso-4}, $s'$ is actually an isomorphism.
	\end{proof}


\bibliography{ref}

\newcommand{\etalchar}[1]{$^{#1}$}
\begin{thebibliography}{ADHM78}

\bibitem[ADHM78]{adhm}
M.F. Atiyah, V.G. Drinfeld, N.J. Hitchin, and Yu.I. Manin.
\newblock Construction of instantons.
\newblock {\em Physics Letters A}, 65(3):185--187, 1978.

\bibitem[AO94]{Hoppe}
V.~Ancona and G.~Ottaviani.
\newblock Stability of special instanton bundles on p2n + 1.
\newblock {\em Transactions of the American Mathematical Society},
  341(2):677--693, 1994.

\bibitem[APR21]{APR}
M.~Altavilla, M.~Petkovic, and F.~Rota.
\newblock Moduli spaces on the kuznetsov component of fano threefolds of index
  2.
\newblock {\em arXiv preprint: 1908.10986}, 2021.

\bibitem[BBF{\etalchar{+}}20]{bayer2020desingularization}
A.~Bayer, S.~Beentjes, S.~Feyzbakhsh, G.~Hein, D.~Martinelli, F.~Rezaee, and
  B.~Schmidt.
\newblock The desingularization of the theta divisor of a cubic threefold as a
  moduli space.
\newblock {\em arXiv preprint arXiv:2011.12240}, 2020.

\bibitem[BBR08]{212cubic}
I.~Biswas, J.~Biswas, and G.V. Ravindra.
\newblock On some moduli spaces of stable vector bundles on cubic and quartic
  threefolds.
\newblock {\em Journal of Pure and Applied Algebra}, 212(10):2298--2306, 2008.

\bibitem[BF08]{BCFacm}
M.~C. Brambilla and D.~Faenzi.
\newblock Moduli spaces of rank-2 acm bundles on prime fano threefolds.
\newblock {\em arXiv preprint: 0806.2265}, 2008.

\bibitem[BF14]{BCDvectorbundlesongenus7}
M.~Chiara Brambilla and D.~Faenzi.
\newblock Vector bundles on fano threefolds of genus 7 and brill–noether
  loci.
\newblock {\em International Journal of Mathematics}, 25(03):1450023, 2014.

\bibitem[BF21]{bolognese2021fullness}
B.~Bolognese and D.~Fiorenza.
\newblock Fullness of exceptional collections via stability conditions--a case
  study: the quadric threefold.
\newblock {\em arXiv preprint arXiv:2103.15205}, 2021.

\bibitem[BLMS17]{BLMS}
A.~Bayer, M.~Lahoz, E.~Macrì, and P.~Stellari.
\newblock Stability conditions on kuznetsov components.
\newblock {\em arXiv preprint: 1703.10839}, 2017.

\bibitem[BM14a]{bayer2014mmp}
A.~Bayer and E.~Macrì.
\newblock Mmp for moduli of sheaves on k3s via wall-crossing: nef and movable
  cones, lagrangian fibrations.
\newblock {\em Inventiones mathematicae}, 198(3):505--590, 2014.

\bibitem[BM14b]{bayer2014projectivity}
A.~Bayer and E.~Macrì.
\newblock Projectivity and birational geometry of bridgeland moduli spaces.
\newblock {\em Journal of the American Mathematical Society}, 27(3):707--752,
  2014.

\bibitem[BMMS12]{BMMS}
M.~Bernardara, E.~Macrì, S.~Mehrotra, and P.~Stellari.
\newblock A categorical invariant for cubic threefolds.
\newblock {\em Advances in Mathematics}, 229:770--803, 2012.

\bibitem[BMS16]{BMS}
A.~Bayer, E.~Macrì, and P.~Stellari.
\newblock The space of stability conditions on abelian threefolds, and on some
  calabi-yau threefolds.
\newblock {\em Invent. math.}, 206:869--933, 2016.

\bibitem[Bri03]{T03}
T.~Bridgeland.
\newblock Stability conditions on triangulated categories.
\newblock {\em Annals of Mathematics}, 166, 01 2003.

\bibitem[DK17]{kronecker}
G.~Dimitrov and L.~Katzarkov.
\newblock Some new categorical invariants.
\newblock {\em arXiv preprint arXiv: 1602.09117}, 2017.

\bibitem[Dru00]{D}
S.~Druel.
\newblock Espace des modules des faisceaux de rang 2 semi-stables de classes de
  chern c1=0, c2=2 et c3=0 sur la cubique de p4.
\newblock {\em International Mathematics Research Notices}, 2000(19):985--1004,
  01 2000.

\bibitem[Fae11]{fae11}
D~Faenzi.
\newblock Even and odd instanton bundles on fano threefolds of picard number 1.
\newblock {\em manuscripta mathematica}, 144, 09 2011.

\bibitem[Har80]{har80}
R.~Hartshorne.
\newblock Stable reflexive sheaves.
\newblock {\em Mathematische Annalen}, 254:121--176, 1980.

\bibitem[HL10]{HL}
D.~Huybrechts and M.~Lehn.
\newblock {\em The Geometry of Moduli Spaces of Sheaves}.
\newblock Cambridge Mathematical Library. Cambridge University Press, 2010.

\bibitem[Hop84]{hop}
H.~J. Hoppe.
\newblock Generischer spaltungstyp und zweite chernklasse stabiler
  vektorraumbündel vom rang 4 auf ip4.
\newblock {\em Mathematische Zeitschrift}, 187:345--360, 1984.

\bibitem[HRS96]{HRS}
D.~Happel, I.~Reiten, and S.~Smalø.
\newblock Tilting in abelian categories and quasitilted algebras.
\newblock {\em Memoirs of the American Mathematical Society}, 575, 03 1996.

\bibitem[IM07]{IM05}
A.~Iliev and L.~Manivel.
\newblock Pfaffian lines and vector bundles on fano threefolds of genus 8.
\newblock {\em Journal of Algebraic Geometry}, 16:499--530, 2007.

\bibitem[JLZ21]{JLZ}
A.~Jacovskis, X.~Lin, and S.~Zhang.
\newblock Hochschild cohomology and categorical torelli for gushel-mukai
  threefolds.
\newblock {\em in preparation}, 2021.

\bibitem[KPS18]{KPS}
A.~Kuznetsov, Y.~Prokhorov, and C.~Shramov.
\newblock Hilbert schemes of lines and conics and automorphism groups of fano
  threefolds.
\newblock {\em Japanese Journal of Mathematics}, 13(1):109–185, Feb 2018.

\bibitem[Kuz07]{kuz06}
A~Kuznetsov.
\newblock Hyperplane sections and derived categories.
\newblock {\em Izvestiya: Mathematics}, 70:447, 10 2007.

\bibitem[Kuz09]{kuz09}
A.~Kuznetsov.
\newblock Derived categories of the fano threefolds v12.
\newblock {\em Proceedings of the Steklov Institute of Mathematics},
  264:110--122, 04 2009.

\bibitem[Kuz12]{kuz12}
A~Kuznetsov.
\newblock Instanton bundles on fano threefolds.
\newblock {\em Central European Journal of Mathematics}, 10:1198--1231, 2012.

\bibitem[Li15]{Li15}
C.~Li.
\newblock Stability conditions on fano threefolds of picard number 1.
\newblock {\em Journal of the European Mathematical Society}, 21, 10 2015.

\bibitem[LMS15]{lahoz2015acm}
M.~Lahoz, E.~Macrì, and P.~Stellari.
\newblock Arithmetically cohen–macaulay bundles on cubic threefolds.
\newblock {\em Algebraic Geometry}, 2:231--269, 2015.

\bibitem[Mac07]{m07}
E.~Macri.
\newblock Stability conditions on curves.
\newblock {\em Mathematical Research Letters}, 14, 06 2007.

\bibitem[Muk92]{mu92}
S.~Mukai.
\newblock Fano 3-folds.
\newblock {\em London Math. Soc. Lecture Note Ser.}, 179, 01 1992.

\bibitem[PR20]{petkovic2020note}
M.~Petkovic and F.~Rota.
\newblock A note on the kuznetsov component of the veronese double cone.
\newblock {\em arXiv preprint arXiv:2007.05555}, 2020.

\bibitem[PR21]{Pertusi2021serreinv}
Laura Pertusi and Ethan Robinett.
\newblock Stability conditions on kuznetsov components of gushel-mukai
  threefolds and serre functor.
\newblock {\em preprint}, 2021.

\bibitem[PY21]{PY20}
L.~Pertusi and S.~Yang.
\newblock {Some Remarks on Fano Three-Folds of Index Two and Stability
  Conditions}.
\newblock {\em International Mathematics Research Notices}, 05 2021.

\bibitem[Qin19]{qin2019moduli}
X.~Qin.
\newblock Moduli space of instanton sheaves on the fano 3-fold v4.
\newblock {\em arXive preprint: 1810.04739}, 2019.

\bibitem[Qin21a]{qin2021bridgeland}
X.~Qin.
\newblock Bridgeland stability of minimal instanton bundles on fano threefolds.
\newblock {\em arXiv preprint arXiv: 2105.14617}, 2021.

\bibitem[Qin21b]{qinv5}
X.~Qin.
\newblock Compactification of the moduli space of minimal instantons on the
  fano 3-fold v5.
\newblock {\em Journal of Pure and Applied Algebra}, 225(3):106526, Mar 2021.

\bibitem[Zha20]{zhang2020bridgeland}
S.~Zhang.
\newblock Bridgeland moduli spaces and kuznetsov's fano threefold conjecture.
\newblock {\em arXiv preprint: 2012.12193}, 2020.

\end{thebibliography}

\bibliographystyle{alpha}

\end{document}